\magnification=1100
\font\bigbf=cmbx10 scaled \magstep2
\font\medbf=cmbx10 scaled \magstep1
\input psfig.sty
\hfuzz=25 pt
\def\sqr#1#2{\vbox{\hrule height .#2pt
\hbox{\vrule width .#2pt height #1pt \kern #1pt
\vrule width .#2pt}\hrule height .#2pt }}
\def\square{\sqr74}
\def\endproof{\hphantom{MM}\hfill\llap{$\square$}\goodbreak}
\def\n{\noindent}

\def\i{\item}

\def\R{I\!\!R}
\def\Z{Z\!\!\!Z}
\def\L{{\bf L}}
\def\T{{\bf T}}
\def\D{{\cal D}}
\def\I{{\cal I}}
\def\forall{\hbox{for all}~}
\def\per{{\rm per}}
\def\sgn{{\rm sign}}
\def\vp{\varphi}
\def\M{{\cal M}}
\def\F{{\cal F}}
\def\O{{\cal O}}
\def\C{{\cal C}}

\def\c{\centerline}
\def\wto{\rightharpoonup}
\def\meas{\hbox{meas}}
\def\v{\vskip 1em}
\def\vs{\vskip 2em}
\def\vsk{\vskip 3em}
\def\ve{\varepsilon}

\def\sign{\,{\rm sign}}
\null
\vs

\c{\bigbf An Optimal Transportation Metric for}
\v
\c{\bigbf Solutions of the Camassa-Holm Equation}
\vs
\c{Alberto Bressan$^{(*)}$ and Massimo Fonte$^{(**)}$}
\vs
\c{(*) Department of Mathematics, Penn State University, University Park, Pa.
16802 U.S.A.}
\v
\c{(**) S.I.S.S.A., Via Beirut 4, Trieste 34014, ITALY}
\vs
\c{Dedicated to Prof.~Joel Smoller in the occasion of his
65-th birthday}
\vsk
\n{\bf Abstract.} In this paper we construct a
global, continuous flow of solutions to the
Camassa-Holm equation on the entire space $H^1$.
Our solutions are conservative, in the sense that
the total energy $\int (u^2+u_x^2)\, dx$ remains a.e.~constant in time.
Our new approach is based on a distance functional $J(u,v)$,
defined in terms of an optimal transportation problem,
which satisfies
${d\over dt} J(u(t), v(t))\leq \kappa\cdot J(u(t),v(t))$ for every
couple of solutions.   Using this new distance functional,
we can construct arbitrary
solutions as the uniform limit of
multi-peakon solutions,
and prove a general uniqueness result.
\vsk
\n{\medbf 1 - Introduction}
\v
The Camassa-Holm equation can be written as a scalar conservation law
with an additional integro-differential term:
$$u_t+(u^2/2)_x+P_x=0\,,\eqno(1.1)$$
where $P$ is defined as a convolution:
$$P\doteq{1\over 2} e^{-|x|} * \left(u^2+{u_x^2\over 2}
\right)\,.\eqno(1.2)$$
For the physical motivations of this equation we refer to
[CH], [CM1], [CM2], [J].  Earlier results on the
existence and uniqueness of solutions can be found in [XZ1], [XZ2].
One can regard (1.1) as an evolution equation on a space of
absolutely continuous functions with derivatives $u_x\in \L^2$.
In the smooth case,
differentiating (1.1) w.r.t.~$x$ one obtains
$$u_{xt}+uu_{xx}+u_x^2-\left( u^2+{u_x^2\over 2}\right)+
P=0\,.\eqno(1.3)$$
Multiplying (1.1) by $u$ and (1.3) by $u_x$ we
obtain the two conservation laws with source term
$$\left(u^2\over 2\right)_t+\left({u^3\over 3}+u\,P\right)_x=u_x \,P\,,
\eqno(1.4)$$
$$\left(u^2_x\over 2\right)_t+\left({uu_x^2\over 2}-{u^3\over 3}\right)_x
=- u_x \,P\,. \eqno(1.5)$$
As a consequence, for regular solutions the total energy
$$E(t)\doteq \int \big[ u^2(t,x)+u_x^2(t,x)\big]\,dx\eqno(1.6)$$
remains constant in time.

As in the case of scalar conservation laws, by
the strong nonlinearity of the equations, solutions
with smooth initial data can lose regularity
in finite time.
For the Camassa-Holm equation (1.1), however,
the uniform bound on $\|u_x\|_{\L^2}$
guarantees that only the $\L^\infty$ norm of the
gradient can blow up, while the solution $u$ itself remains
H\"older continuous at all times.

In order to construct global in time solutions, two
main approaches have recently been introduced.
On one hand, one can add a small diffusion term in the
right hand side of (1.1), and recover solutions of the
original equations as a vanishing viscosity limit [CHK1, CHK2].
An alternative technique, developed in [BC2], relies on
a new set of independent and dependent variables, specifically
designed with the aim of
``resolving'' all singularities.
In terms of these new variables, the solution
to the Cauchy problem becomes regular for all times, and can
be obtained as the unique fixed point of a contractive transformation.

In the present paper, we implement yet another approach
to the Camassa-Holm equation.
As a starting point we consider all multi-peakon solutions, of the form
$$u(t,x)=\sum_{i=1}^N p_i(t) e^{-|x-q_i(t)|}\,.\eqno(1.7)$$
These are obtained by solving the system of O.D.E's
$$\left\{
\eqalign{ \dot q_i&=\sum_j p_j\,e^{-|q_i-q_j|}\,,\cr
\dot p_i&=\sum_{j\not= i} p_i p_j\,\sgn(q_i-q_j)\,e^{-|q_i-q_j|}\,.
\cr}\right.\eqno(1.8)$$
It is well known that this can
be written in hamiltonian form:
$$\left\{
\eqalign{ \dot q_i&~=~~{\partial\over\partial p_i} H(p,q)\,,
\cr
\dot p_i&~=~-{\partial\over\partial q_i} H(p,q)\,,
\cr}\right.\qquad\qquad H(p,q)\doteq {1\over 2} \sum_{i,j} p_ip_j
e^{-|q_i-q_j|}\,.
$$
If all the coefficients $p_i$ are initially positive, then
they remain positive and bounded for all times.
The solution $u=u(t,x)$ is thus uniformly Lipschitz continuous.
We stress, however, that
here we are not making any assumption about the signs of the
$p_i$.   In a typical situation,
two peakons can
cross at a finite time $\tau$. As $t\to\tau-$
their strengths $p_i,p_j$ and positions $q_i,q_j$
will satisfy
$$p_i(t)\to +\infty\,,\qquad p_j(t)\to -\infty\,,
\qquad p_i(t)+p_j(t)\to \bar p\,,
\eqno(1.9)$$
$$q_i(t)\to \bar q\,,
\quad q_j(t)\to \bar q\,,\qquad q_i(t)<q_j(t)
~~\hbox{for} ~~t<\tau\,,\eqno(1.10)$$
for some $\bar p,\bar q\in\R$.
Moreover, $\big\|u_x(t)\big\|_{\L^\infty}\to\infty$.
In this case, we will show that
there exists
a unique way to extend the multi-peakon solution beyond
the interaction time, so that the total energy is conserved.

Having constructed a set of ``multi-peakon solutions'',
our main goal is to show that these solutions form a
continuous semigroup, whose domain is dense in the space
$H^1(\R)$.  Taking the unique continuous extension,
we thus obtain a continuous semigroup of solutions of
(1.1), defined on the entire space $H^1$.

The particular metric, used to derive the Lipschitz continuity of
the semigroup, is the most novel feature of our paper.
Indeed, one easily checks that the flow map
$\Phi_t: u(0)\mapsto u(t)$ cannot be continuous as
a map from $H^1$ into itself, or from $\L^2$ into itself.
Distances defined in terms of convex norms perform well
in connection with linear problems, but occasionally fail
when nonlinear features become dominant.
In the present setting, we construct a new distance
$J(u,v)$ between functions $u,v\in H^1$, defined
by a problem of optimal transportation. Roughly speaking,
$J(u,v)$ will be the minimum cost in order to transport the
mass distribution with density $1+u_x^2$ located
on the graph of $u$ onto the
mass distribution with density $1+v_x^2$ located
on the graph of $v$.  See Section 3 for details.
With this definition of distance, our main result shows that
$$\left|{d\over dt} J\big(u(t),\,v(t)\big)\right|~\leq~ C \cdot
J\big(u(t),\,v(t)\big)\eqno(1.11)$$
for some constant $C$ and any couple of multi-peakon solutions $u,v$.
Moreover, $J(u_n,u)\to 0$ implies the uniform convergence
$\|u_n- u\|_{\L^\infty}\to 0$.
The distance functional $J$ thus provides the ideal tool to measure
continuous dependence on the initial data for solutions
to the Camassa-Holm equation.
Earlier applications
of distances defined in terms of optimal transportation
problems can be found in the monograph [V].

The issue of uniqueness of solutions must here be discussed in
greater detail.  For a multi-peakon solution, as long as all
coefficients $p_i$ remain bounded, the solution to
the system of ODE's  (1.8) is clearly unique.   For each time $t$, call
$\mu_t$ the measure having density
$u^2(t)+u_x^2(t)$ w.r.t.~Lebesgue measure.
Consider a time $\tau$
where a positive and a negative peakon collide, according
to (1.9)-(1.10).
As $t\to \tau- $, we have the weak convergence $\mu_t\wto
\mu_\tau$ for some positive measure $\mu_\tau$ which typically
contains a Dirac mass at the point $\bar q$.
By energy conservation, we thus have
$$\int\big[ u^2(\tau,x)+u_x^2(\tau, x)\big]\,dx +\mu_\tau\big(\{\bar q\}
\big)=\lim_{t\to\tau-}\int\big[ u^2(\tau,x)+u_x^2(\tau, x)\big]\,dx
= E(\tau-)\,.$$
There are now two natural ways to prolong the multi-peakon solution
beyond time $\tau$:   a conservative solution, such that
$$E(t)=\int\big[ u^2(t,x)+u_x^2(t, x)\big]\,dx = E(\tau -)\qquad\qquad
t>\tau\,,$$
or a dissipative solution, where all the energy concentrated
at the point $\bar q$ is lost.
In this case
$$E(t)=\int\big[ u^2(t,x)+u_x^2(t, x)\big]\,dx = E(\tau -)-
\mu_\tau\big(\{\bar q\}\big)\qquad\qquad
t>\tau\,.$$
For $t>\tau$,
the dissipative solution is obtained by simply replacing the two peakons
$p_i,p_j$ with one single peakon of strength $\bar p$, located at
$x=\bar q$.  On the other hand, as we will show in Section 2,
the conservative solution contains two peakons emerging from
the point $\bar q$.  As $t\to \tau +$, their strengths and positions
satisfy again (1.9), while (1.10) is replaced by
$$q_i(t)\to \bar q\,,
\quad q_j(t)\to \bar q\,,\qquad q_i(t)>q_j(t)
~~\hbox{for} ~~t>\tau\,.\eqno(1.12)$$
The vanishing viscosity approach in [CHK1, CHK2] singles out
the dissipative solutions.  These can also be characterized
by the Oleinik type estimate
$$u_x(t,x)\leq C(1+t^{-1})\,,$$
valid for $t>0$ at a.e.~$x\in\R$.
On the other hand, the coordinate transformation approach
in [BC2] and the present one, based on optimal transport metrics,
appear to be well suited for the study of
both conservative and dissipative solutions.

In the present paper we focus on conservative solutions
to the C-H equation in the spatially periodic case.
The assumption of spatial periodicity allows us to concentrate
on the heart of the matter, i.e.~the uniqueness and stability
of solutions beyond the time of singularity formation.
It will spare us some technicalities, such as the analysis of the
tail decay of $u, u_x$ as $x\to\pm\infty$.
The construction of an optimal transportation metric
in connection with dissipative solutions to the C-H equation
is significantly different, and
will be carried out in a forthcoming paper. The main ingredients
can already be found in the paper [BC1], devoted to dissipative solutions
of the Hunter-Saxton equation.

As initial data, we take
$$u(0,x)=\bar u(x)\,,\eqno(1.13)$$
with $\bar u$ in the space $H^1_\per $ of
of periodic, absolutely
continuous
functions $u$ with derivative $u_x\in\L^2_\per $.
To fix the ideas, we assume that the period is $1$,
so that
$$u(x+1)=u(x)\qquad\qquad ~x\in\R\,.$$
On $H^1_\per $ we shall use the norm
$$\big\|u\|_{H^1_\per }\doteq \left(\int_0^1\big|u(x)\big|^2\,dx
+\int_0^1\big|u_x(x)\big|^2\,dx\right)^{1/2}.$$
Our main results can be stated as follows.
\v
\n{\bf Theorem 1.} {\it For each initial
data $\bar u\in H^1_\per $, there exists
a solution $u(\cdot)$ of the Cauchy
problem (1.1), (1.13). Namely, the map
$t\mapsto u(t)$ is Lipschitz continuous
from $\R$ into $\L^2_\per $, satisfies (1.13) at time $t=0$,
and the identity
$${d\over dt} u = -uu_x -P_x\eqno(1.14)$$
is satisfied as an equality between elements in $\L^2_\per $
at a.e.~time $t\in\R$. This same map $t\mapsto u(t)$ is continuously
differentiable from $\R$ into $\L^p_\per$ and satisfies (1.14)
at a.e.~time $t\in\R$, for all $p\in [1,2[\,$.
The above solution is conservative in the sense that, for a.e.~$t\in\R$,
$$E(t)=\int_0^1 \big[u^2(t,x)+u_x^2(t,x)\big]\,dx=E^{\bar u}\doteq
\int_0^1 \big[\bar u^2(x)+\bar u_x^2(x)\big]\,dx\,.\eqno(1.15)$$
}
\v
\n{\bf Theorem 2.} {\it Conservative solutions to (1.1) can be constructed
so that they constitute a continuous flow $\Phi$. Namely, there exists
a distance functional $J$ on $H^1_\per $ such that
$${1\over C}\cdot\|u-v\|_{\L^1_\per}\leq J(u,v)\leq C\cdot
\|u-v\|_{H^1_\per }
\eqno(1.16)$$
for all $u,v\in H^1_\per$ and some constant $C$ uniformly valid on
bounded sets of $H^1_\per\,$.
Moreover, for any two solutions
$u(t)=\Phi_t\bar u$, $v(t)=\Phi_t\bar v$ of (1.1),
the map $t\mapsto J\big(u(t),\,v(t)\big)$
satisfies
$$J\big(u(t),\,\bar u\big)\leq C_1\cdot |t|\,,\eqno(1.17)$$
$$J\big(u(t),\,v(t)\big)~\leq~ J(\bar u,\bar v)\cdot e^{C_2 |t|}
\eqno(1.18)$$
for a.e. $t\in\R\,$ and constants $C_1,C_2$, uniformly valid as
$u,v$ range on bounded sets of $H^1_{per}\,$.
}
\v
Somewhat surprisingly, all the properties stated in Theorem 1
are still not strong enough
to single out a unique solution.  To achieve uniqueness, an
additional condition is needed.
\v
\n{\bf Theorem 3.} {\it Conservative solutions $t\mapsto u(t) $ of (1.1)
can be constructed with the following additional property:
\v
\n For each $t\in\R$, call $\mu_t$ the absolutely continuous
measure having density $u^2+u_x^2$ w.r.t.~Lebesgue measure.
Then, by possibly redefining $\mu_t$ on a set of times of measure zero,
the map $t\mapsto\mu_t$ is continuous w.r.t.~the topology
of weak convergence of measures. It provides a measure-valued
solution to the conservation law
$$w_t+(uw)_x=(u^3-2uP^u)_x\,.\eqno(1.19)$$
\v
The solution of the Cauchy problem (1.1), (1.13) satisfying the properties
stated in Theorem 1 and this additional condition is unique.}
\v
The plan of the paper is as follows.  In Section 2 we derive some
elementary properties
of multi-peakon solutions and show that any initial data can be
approximated in $H^1_\per$ by a finite sum of peakons.
In Section 3 we introduce our distance functional $J(u,v)$
and study its
relations with other distances defined by Sobolev norms.
The continuity of the flow (1.1), together with the key
estimates (1.17)-(1.18)
are then proved in the following two sections.
The proofs of Theorems 1 and 2 are completed in
Section 6.  The uniqueness result stated in Theorem 3 is proved
in Section 7.  As a corollary, we also show that in a multi-peakon
solution the only possible interactions
involve exactly two peakons: one positive and one negative.
In particular, no triple interactions can ever occur.
\vsk

\n{\medbf 2 - Multipeakon solutions}
\v
By a periodic peakon we mean a function of the form
$$u(x)=p\, \chi(x-q)\,,\qquad\qquad
\chi(x)\doteq\sum_{n\in\Z} e^{-|x-n|}\,.\eqno(2.1)$$
Observe that the periodic function $\chi$ satisfies
$$
\eqalign{
&\chi(-x)=\chi(x)=\chi(x+1)\qquad\qquad x\in\R\,,
\cr
&\chi(x)
= {e^x+e^{1-x}\over e-1}\qquad\qquad\quad\qquad x\in [0,1]\,.
}\eqno(2.2)$$

We begin this section by observing that any periodic initial data
can be approximated by a periodic multi-peakon.
\v
\n{\bf Lemma 1.} {\it Let $f\in H^1_\per$.
Then for any $\ve>0$
there exists periodic multi-peakon $g$,
of the form
$$g(x)=\sum_{i=1}^{N} p_i\sum_{n\in\Z} e^{-|x-q_i-n|}=
\sum_{i=1}^N p_i\,\chi(x-q_i)\eqno(2.3)$$
such that}
$$\|f-g\|_{H^1_\per}<\ve\,.$$
\v
\n{\bf Proof.} By taking a suitable mollification, we can
approximate $f$ with a periodic function $\tilde f\in\C^\infty$, so that
$$\|f-\tilde f\|_{H^1_\per}<\ve/2\,.\eqno(2.4)$$
Next, we observe that
$${1\over 2}\left( e^{-|x|}-{\partial^2\over\partial x^2} e^{-|x|}
\right)=\delta_0\,,$$
where $\delta_0$ denotes the Dirac distribution concentrating
a unit mass at the origin.
We can thus write $\tilde f$ as a convolution:
$$
\tilde f=\delta_0*\tilde f ={1\over 2}
\left( e^{-|x|}-{\partial^2\over\partial x^2} e^{-|x|}
\right)*\tilde f = e^{-|x|} * \left({\tilde f-\tilde f''\over 2}\right)\,,
$$
$$\tilde f(x)=\int_0^1\chi(x-y)\cdot {\tilde f(y)-\tilde f''(y)\over 2}\,dy\,.$$
The above integral
can now be approximated with a Riemann sum
$$
g(x)=\sum_{i=1}^N p_i\,\chi(x-q_i)\,,
\qquad\qquad p_i=\int_{(i-1)/N}^{i/N}{\tilde f(y)-\tilde f''(y)\over 2}\,dy\,.$$
Choosing $N$ sufficiently large we obtain
$\|\tilde f-g\|_{H^1_{per}}<\ve/2$. Together with (2.4)
this yields the result.
\endproof
\vs
Next, we show how to construct a unique conservative solution,
for multi-peakon initial data.  As long as the locations
$q_i$ of the peakons remain distinct, this can be obtained by
solving the Hamiltonian system of O.D.E's (1.8).
However, at a time $\tau$ where two or more
peakons interact, the corresponding strengths
$p_i$ become unbounded.  A suitable transformation of
variables is needed, in order to resolve the singularity
and uniquely extend the solution beyond the interaction time.
\v
\n{\bf Lemma 2.} {\it Let $\bar u$ be any periodic, multi-peakon
initial data. Then the Cauchy problem (1.1), (1.13)
has a global, conservative multi-peakon solution defined for all $t\in\R$.
The set $\I$ of times where two or more peakons interact is at most
countable.
Moreover, for all $t\notin\I$, the energy conservation (1.15) holds.
}
\v
\n {\bf Proof.} The solution can be uniquely constructed by
solving the hamiltonian system (1.8), up to the first time
$\tau$ where two or more peakons interact.
We now show that there exists a unique way to prolong the solution
for $t>\tau$, in terms of two outgoing peakons.
To fix the ideas, call
$$\bar q=\lim_{t\to\tau-} q_i(t)\qquad\qquad i=1,\ldots,k\,,$$
the place
where the interaction occurs,
and let $p_1(t),\ldots, p_k(t)$ be the strengths of the
interacting peakons.
Later in this paper we will show that only the case $k=2$
can actually occur, but at this stage we need to consider
the more general case.
We observe that the strengths $p_{k+1},\ldots,p_N$ of the peakons
not involved in the
interaction remain continuous at time $\tau$.   Moreover,
by (1.8)
there exists the limit
$$\bar p=\lim_{t\to \tau-} \sum_{i=1}^k p_i(t)\,.$$
We can thus write
$$u(\tau,x)~=~\lim_{t\to\tau-}
\sum_{i=1}^N p_i(t)\,e^{-|x-q_i(t)|}~=~\bar p \,
e^{-|x-\bar q|}+\sum_{i=k+1}^N
p_i(\tau)\,e^{-|x-q_i(\tau)|}\,.$$
For $t>\tau$, we shall prolong the solution with two
peakons emerging from the point $\bar q$.  The strength
of these two peakons will be uniquely determined
by the requirement of energy conservation (1.15).

Call $\xi^-(t)$, $\xi^+(t)$ respectively
the position of the smallest and largest characteristic
curves passing through the point $(\tau,\bar q)$, namely
$$\eqalign{\xi^-(t)&\doteq \min\Big\{ \xi(t)\,;~~\xi(\tau)=\bar q\,,
~~\dot\xi(s)=
u\big(s,\xi(s)\big)\qquad\forall s\in [\tau-h,\tau+h]\Big\}\,,\cr
\xi^+(t)&\doteq \max\Big\{ \xi(t)\,;~~\xi(\tau)=\bar q\,,
~~\dot\xi(s)=
u\big(s,\xi(s)\big)\qquad\forall s\in [\tau-h,\tau+h]\Big\}\,.\cr}
\eqno(2.5)$$
Moreover, define
$$e_{(\tau,\bar q)}\doteq\lim_{t\to \tau-}
\int_{\xi^-(t)}^{\xi^+(t)} u_x^2(t,x)\,dx\,.\eqno(2.6)$$
The existence of this limit follows from the balance law
(1.5).
This describes how much energy is concentrated at the interaction point.

For $t>\tau$ the solution will contain the peakons
$p_{k+1},\ldots,p_N$, located at $q_{k+1},\ldots, q_N$, together with
the two outgoing peakons $p_1,p_2$,  located at $q_1<q_2$.
The behavior of $p_i,q_i$ for $i\in\{k+1,\ldots, N\}$ is
still described by a system of O.D.E's as in (1.8).
However, to describe the evolution of $p_1,p_2,q_1,q_2$ one has to
use a different set of variables, resolving the singularity
occurring at $(\tau,\bar q)$.
As $t\to \tau+$ we expect
(1.9), (1.12) to hold.
To devise a suitable set of rescaled variables,
we observe that, by (1.3),
$${d\over dt}\,u_x\big(t,\xi(t)\big)= -{1\over 2}
\,u^2_x\big(t,\xi(t)\big) + [u^2-P]\eqno(2.7)$$
along any characteristic curve $t\mapsto \xi(t)$
emerging from the point $\bar q$.
Since $u,P$ remain uniformly bounded, one has
$$u_x(t,x)\approx {2\over t-\tau}\qquad\qquad t>\tau\,,~~x\in \big[
q_1(t)\,,~q_2(t)\big]\,.$$
The total amount of energy concentrated in the interval between
the two peakons is given by
$$\eqalign{\int_{q_1(t)}^{q_2(t)} \big[ u^2(t,x)+u_x^2(t,x)\big]\,dx
&\approx \left({u(q_2)- u(q_1)\over q_2-q_1}\right)^2\cdot (q_2-q_1)
\approx {\Big[ (p_2-p_1)\big(1-e^{-|q_2-q_1|}\big)\Big]^2\over q_2-q_1}
\cr
&\approx (p_2-p_1)^2(q_2-q_1)
\approx e_{(\tau,\bar q)}\,.\cr}$$
The previous heuristic analysis suggests that, in order to
resolve the singularities, we should work with the variables
$$z=p_1+p_2\,,\qquad w=2\arctan (p_2-p_1)\,,\qquad
\eta=q_2+q_1\,,\qquad \zeta =(p_2-p_1)^2(q_2-q_1)\,,$$
together with $p_{k+1}\,,\ldots,\,p_N\,,~~~q_{k+1}\,,\ldots, \, q_N$.
{}
To simplify the following calculations we here
assume $0<q_1<q_2<q_{k+1}<...<q_N<1$,
which is not restrictive.

Let $\chi$ defined in (2.1) and $\tilde \chi(x)\doteq {-e^x+e^{1-x}\over e-1}$, $x\in [0,1]$.
From the original system of equations (1.8) it follows
$$
\eqalign{
&\dot z=
\cosh\left({\zeta\over 2 \tan^2{w\over 2}}\right)
z\cos^2{w\over 2}\sum_{j=k+1}^N p_j \chi\left(q_j-{\eta\over 2}\right)
-{\sinh \left({\zeta\over 2 \tan^2{w\over 2}}\right)\over
{1\over \tan {w\over 2}}}\sum_{j=k+1}^N p_j \tilde\chi\left(q_j-{\eta\over 2}
\right)
\cr
&\dot w=\left(z^2\cos^2 {w\over 2}-\sin^2{w\over 2}\right)
\chi\left({\zeta\over \tan^2{w\over 2}}\right) +
2
\cosh \left({\zeta\over 2 \tan^2{w\over 2}}\right)z
\sum_{j=k+1}^N p_j \chi\left(q_j-{\eta\over 2}\right)
\cr
&\qquad
+
2
\sinh \left({\zeta\over 2 \tan^2{w\over 2}}\right)\sin w
\sum_{j=k+1}^N p_j \tilde\chi\left(q_j-{\eta\over 2}\right)
\cr
&\dot \eta=
z\left[\chi(0)+\chi\left({\zeta\over \tan^2{w\over 2}}\right)\right]+
\cosh\left({\zeta\over 2\tan^2{w\over 2}}\right)\sum_{j=k+1}^N p_j \chi\left(q_j-{\eta\over 2}\right)
\cr
&\dot \zeta=
{\chi(0)-\chi\left({\zeta\over \tan^2{w\over 2}}\right)\over {1\over \tan^2{w\over 2}}}\zeta+
\chi\left({\zeta\over \tan^2{w\over 2}}\right){z^2\zeta\over \tan{w\over 2}}
-{\sinh\left(\zeta\over 2 \tan^2{w\over 2}\right)\over {1\over \tan^2{w\over 2}}}
\sum_{j=k+1}^N p_j\tilde\chi\left(q_j-{\eta\over 2}\right)
\cr
&\qquad+2{\zeta\over \tan{w\over 2}}
\left[
\cosh\left(\zeta\over 2\tan^2{w\over 2}\right){z \over \tan{w\over 2}}
\sum_{j=k+1}^N p_j\chi\left(q_j-{\eta\over 2}\right)
-
\sinh\left(\zeta\over 2\tan^2{w\over 2}\right)
\sum_{j=k+1}^N p_j\tilde\chi\left(q_j-{\eta\over 2}\right)
\right]
\cr
&\dot p_i =p_i
\left[
\cosh\left({\zeta\over 2\tan^2{w\over 2}}\right)z
\chi\left(q_i-{\eta\over 2}\right)
+
\sinh\left({\zeta\over 2\tan^2{w\over 2}}\right)\tan{w\over 2}\tilde
\chi\left(q_i-{\eta\over 2}\right)
\right]
\cr
&\qquad+p_i\sum_{j=k+1}^N p_j \sign(q_i- q_j)\chi\left({|q_i-q_j|}\right)
\cr
&\dot q_i=
\cosh\left({\zeta\over 2\tan^2{w\over 2}}\right)z
\chi\left(q_i-{\eta\over 2}\right)
+
\sinh\left({\zeta\over 2\tan^2{w\over 2}}\right)\tan{w\over 2}\tilde \chi
\left(q_i-{\eta\over 2}\right)
+\sum_{j=k+1}^N p_j\chi\left({|q_i-q_j|}\right)
}
$$
with initial data
$$z(\tau)= \bar p\,,\qquad w(\tau)=\pi\,,\qquad
\eta(\tau)= 2\bar q\,,\qquad \zeta(\tau)= e_{(\tau,\bar q)}\,,$$
$$p_i(\tau)=\lim_{t\to \tau-}p_i(t)\,,
\qquad\qquad q_i(\tau)=\lim_{t\to \tau-}q_i(t)\qquad\qquad i=k+1,\ldots,N\,.$$
For the above system of O.D.E's, a direct inspection
reveals that the right hand side can be
extended by continuity also at the value $w=\pi$, because
all singularities are removable.
This continuous extension is actually smooth,
in a neighborhood of the initial data. Therefore, our Cauchy
problem has a unique local solution.
This provides a multi-peakon solution defined on some interval of the
form $[\tau,\,\tau'[\,$, up to the next interaction time.

The case where two or more groups of peakons interact
exactly at the same time $\tau$, but at different locations within
the interval $[0,1]$,
can be treated in exactly
the same way.
Since the total number of peakons (on a unit interval in the
$x$-variable) does not increase, it is clear that the number of
interaction times is at most countable.  The solution can thus be extended
to all times $t>0$, conserving its total energy.
\endproof

\vsk
\n{\medbf 3 - A distance functional}
\v
In this section we shall
construct a functional $J(u,v)$ which controls the distance
between two solutions of the equation (1.1).
All functions and measures on $\R$ are assumed to be periodic with
period 1.
Let $\T$ be the unit circle, so that $\T=[0,2\pi]$ with the
endpoints $0$ and $2\pi$ identified.
The distance $|\theta-\theta'|_*$ between two points $\theta,\theta'\in\T$
is defined as the smaller between the lengths of the two arcs connecting
$\theta$ with $\theta'$ (one clockwise, the other counterclockwise).
We now consider the product space
$$X\doteq \R\times\R\times \T$$
with distance
$$d^\diamondsuit\Big( (x,u,w), ~(\tilde x, \tilde u,\tilde w)\Big)
\doteq \Big(|x-\tilde x|+|u-\tilde u|+ |w-\tilde
w|_*\Big)\wedge 1 \,,\eqno(3.1)$$
where $a\wedge b\doteq \min\{a,b\}$.
Let $\M(X)$ be the space of all Radon measures on $X$ which are 1-periodic
w.r.t.~the $x$-variable. To
each 1-periodic function $u\in H^1_\per $
we now associate the positive measure
$\sigma^u\in\M(X)$ defined as
$$\sigma^{u}(A)\doteq \int_{\big\{ x\in\R\,:~(x,\,u(x),\, 2\arctan
\,u_x(x)\,)\in A\big\}} \big(1+u_x^2(x)\big)\, dx
\eqno(3.2)$$
for every Borel
set $A\subseteq\R^2\times \T\,$.
Notice that the total mass of $\sigma^{(u,\mu)}$ over one period is
$$\sigma^u\big([0,1]\times\R\times\T)=
1+\int_0^1 u_x^2(x)\,dx
\,.$$

On this family of positive, 1-periodic Radon measures, we now introduce
a kind of Kantorovich distance,
related to an optimal transportation problem.
Given the two measures $\sigma^u$ and
$\sigma^{\tilde u}$,
their distance $J(u,\tilde u)$ is defined as follows.
\v
Call $\F$ the family of all
strictly increasing absolutely continuous
maps $\psi:\R\mapsto\R$ which have an absolutely continuous
inverse and satisfy
$$\psi(x+n)=n+\psi(x)\qquad\qquad \hbox{for every}~~n\in \Z\,.\eqno(3.3)$$
For a given $\psi\in\F$, we define the 1-periodic, measurable functions
$\phi_1,\phi_2:\R\mapsto [0,1]$ by setting
$$\eqalign{\phi_1(x)&\doteq
\sup\,\bigg\{\theta\in [0,1]\,;~~ \theta\cdot
\Big( 1+u_x^2(x)\Big)\leq \Big(
1+ \tilde u_x^2
\big(\psi(x)\big)\Big)\,\psi'(x)\bigg\}\,,\cr
\phi_2(x)&\doteq
\sup\,\bigg\{\theta\in [0,1]\,;~~
1+u_x^2(x)\geq \theta\cdot \Big(
1+ \tilde u_x^2
\big(\psi(x)\big)\Big)\,\psi'(x)\bigg\}\,.\cr
}\eqno(3.4)$$
Observe that the above definitions imply
$\max\big\{ \phi_1(x),\,\phi_2(x)\big\}=1$ together with
$$\phi_1(x) \,\Big( 1+u_x^2(x)\Big)=\phi_2\big(\psi(x)\big)\,\Big(
1+\tilde u_x^2
\big(\psi(x)\big)\Big)\,\psi'(x)\eqno(3.5)$$
for a.e. $x\in\R$.
We now define
$$\eqalign{J^\psi(u,\tilde u)\doteq &
\int_0^1 d^\diamondsuit\Big( \big(x,\,u(x),\, 2\arctan u_x(x)\big)\,,~
\big(\psi(x),\,\tilde u(\psi(x)),\,2\arctan \tilde u_x(\psi(x))\Big) \cr
&\qquad\qquad\qquad\qquad\qquad\qquad \cdot
\phi_1(x)\,\big(1+u_x^2(x)\big)\,dx\cr
&+\int_0^1\Big| \big(1+u_x^2(x)\big)-\big(1+\tilde u_x^2(\psi(x))\big)
\,\psi'(x)\Big|\,dx\,.\cr}\eqno(3.6)$$
Of course, the integral is always computed over one period.
Observe that $x\mapsto \psi(x)$ can be
regarded as a {\bf transportation plan}, in order to transport the
measure $\sigma^u$ onto the measure $\sigma^{\tilde u}$.
Since these two positive
measures need not have the same total mass, we allow the presence of
some excess mass, not transferred
from one place to the other. The penalty for this
excess mass is given by the second integral in (3.6).
The factor $\phi_1\leq 1$ in the first integral indicates
the percentage of the mass which is actually transported.
Integrating (3.5)
over one period, we find
$$\int_0^1\phi_1(x) \big(1+u_x^2(x)\big)\,dx=
\int_0^1 \phi_2(y) \big(1+\tilde u_x^2(y)\big)\,dy\,.$$
We can thus transport the
measure $\phi_1 \,\sigma^u$ onto $\phi_2\,
\sigma^{\tilde u}$ by a map
$\Psi:\big(x,\, u(x)\, \arctan u_x(x)\big)\mapsto \big(y,
\,\tilde u(y),\,\arctan \tilde u_x(y)\big)$, with $y=\psi(x)$.
The associated
cost is given by the first integral in (3.6). Notice that in this
case the measure $\phi_2\, \sigma^{\tilde u}$
is obtained as the push-forward
of the measure $\phi_1 \,\sigma^u$. We recall that the {\bf
push-forward} of a measure $\sigma$ by a mapping $\Psi$ is defined as
$(\Psi\#\sigma)(A)\doteq \sigma(\Psi^{-1}(A))$ for every measurable set
$A$.  Here $\Psi^{-1}(A)\doteq \big\{ z\,;~\Psi(z)\in A\big\}$.

\v Our distance functional $J$ is now obtained by optimizing over all
transportation plans, namely
$$J(u,\tilde u)\doteq \inf_{\psi\in\F} J^\psi
(u,\tilde u)\,.\eqno(3.7)$$
\v
To check
that (3.7) actually defines a distance, let $u,v,w\in H^1(\R)$ be
given.
\v
\n{\bf 1.} Choosing $\psi(x)=x$, so that
$\phi_1(x)=\phi_2(x)=1$, we immediately see that $J(u,u)=0$. Moreover,
if $J(u,\tilde u)=0$, then by the definition of $d^\diamondsuit$
we have $\tilde u= u$.
\v
\n{\bf 2.} Given $\psi\in\F$,
define $\tilde\psi=\psi^{-1}$, so that $\tilde\phi_1=\phi_2$,
$\tilde\phi_2=\phi_1$. This yields
$$J^{\tilde \psi}(\tilde u, u)=
 J^\psi(u,\tilde u)\,.$$
Hence $J(\tilde u,u)=J(u,\tilde u)$.
\v
\n{\bf 3.} Finally, to prove the
triangle inequality, let $\psi^\flat,\psi^\sharp:\R\mapsto\R$ be
two increasing diffeomorphisms satisfying (3.3), and let
$\phi_1^\flat,\phi_1^\sharp,\phi_2^\flat,\phi_2^\sharp:\R\mapsto
[0,1]$ be the corresponding functions, defined as in (3.4).
We now consider the
composition $\psi\doteq \psi^\sharp\circ\psi^\flat$ and define
the functions $\phi_1,\phi_2$ according to (3.4).
Observing that
$$\phi_1(x)\geq \phi_1^\flat(x)\cdot\phi_1^\sharp\big(\psi^\flat(x)\big)
\,,$$
$$\phi_2\big(\psi(x)\big)=\phi_2\Big(\psi^\sharp\big(\psi^\flat(x)\big)\Big)
\geq
\phi_2^\flat\big(\psi^\flat(x)\big)\cdot\phi_2^\sharp
\Big(\psi^\sharp\big(\psi^\flat(x)\big)\Big)
\,,$$
and recalling that the distance $d^\diamondsuit$ at (3.1) is always $\leq 1$,
we conclude
$$J^\psi(u,w)\leq
J^{\psi^\flat}(u,v)+J^{\psi^\sharp}(v,w)\,.
$$
This implies the triangle inequality $J(u,v)+J(v,w)\geq
J(u,w)$.
\endproof

In the remainder of this section we study the relations between our
distance functional $J$ and the distances determined by various norms.
\v
\n{\bf Lemma 3.} {\it For any $u,v\in H^1_\per$ one has
$${1\over C}\cdot
\|u-v\|_{\L^1_\per}\leq J(u,v)\leq C\cdot\|u-v\|_{H^1_\per}\,,
\eqno(3.8)$$
with a constant $C$ uniformly valid on bounded subsets of $H^1_\per$.
}
\v
\n{\bf Proof.}
We shall use the elementary bound
$$\big| \arctan a-\arctan b\big|\cdot a^2\leq
4\pi\big(|a|+|b|\big)\,|a-b|\,,\eqno(3.9)$$
valid for all $a,b\in\R$.
In connection with the identity mapping $\psi(x)=x$
we now compute
$$\eqalign{ J^\psi(u,v)&\leq \int_0^1 \Big\{
\big| u(x)-v(x)\big|+2 \big|\arctan u_x-\arctan v_x\big|
\Big\}\, (1+u_x^2)\,dx
+\int_0^1 \big|u_x^2-v_x^2|\,dx\cr
&\leq \|u-v\|_{\L^\infty}\,\|1+u_x^2\|_{\L^1}+ (8\pi+1)\int_0^1
|u_x+v_x|\, |u_x-v_x|\,dx\cr
&\leq (8\pi + 3) \,\big(1+ \|u\|_{H^1}+\|v\|_{H^1}\big)\cdot
\|u-v\|_{H^1}\,,\cr}$$
proving the second inequality in (3.8).

\midinsert
\vskip 10pt
\centerline{\hbox{\psfig{figure=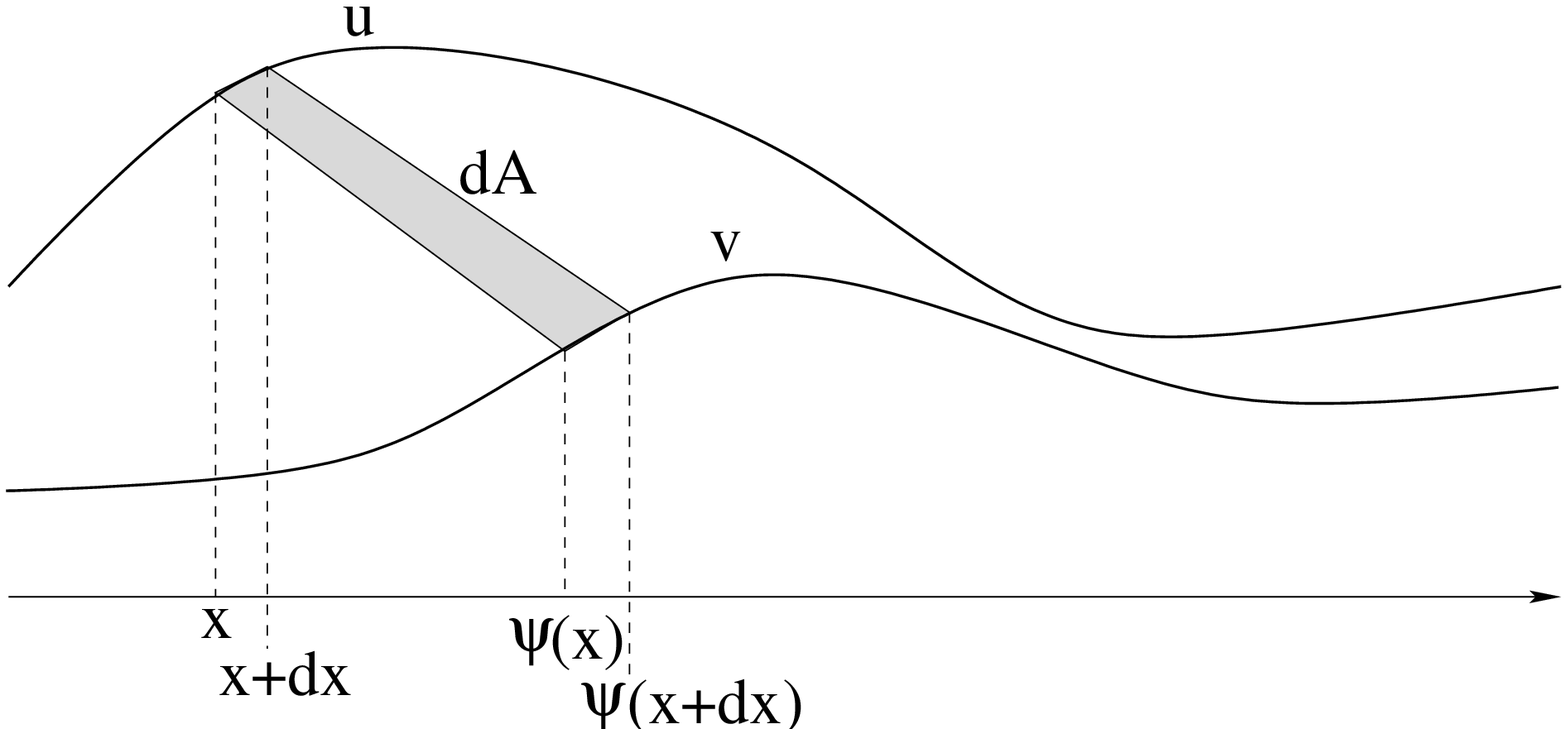,width=10cm}}}
\centerline{\hbox{figure 1}}
\vskip 10pt
\endinsert

To achieve the first inequality, choose any $\psi\in\F$.
For $x\in [0,1]$, call $\gamma^x$ the segment joining the
point $P^x=\big(x,u(x)\big)$ with $Q^x=\big(\psi(x), v(\psi(x))\big)$.
Clearly, the union of all these segments covers the
region between the graphs of $u$ and $v$.
Moving the base point from $x$ to $x+dx$, the corresponding segments
sweep an infinitesimal area $dA$ estimated by (fig.~1)
$$\eqalign{|dA|&\leq |P^x-Q^x|\cdot \big(|dP^x|+|dQ^x|\big)
\cr
&\leq
\Big(\big|x-\psi(x)\big|^2+\big|u(x)-v(\psi(x))\big|^2+\Big)^{1/2}
\cdot
\Big[(1+u_x^2)^{1/2}dx+(1+v_x^2)^{1/2}\psi'(x)\,dx\Big]\,.\cr}
$$
Integrating over one period we obtain
$$\eqalign{ \int_0^1&\big|u(x)-v(x)\big|\,dx\cr
&\leq
\int_0^1\Big(\big|x-\psi_{(t)}(x)\big|+
\big|u(x)-v(\psi_{(t)}(x))\big|\Big)
\cdot
\Big[\big(1+u_x^2(x)\big)+\big(1+v_x^2(\psi(x))\big)^{1/2}
\psi'_{(t)}(x)\Big]\,dx
\cr
&\leq \big(2+\|u\|_{H^1}+\|v\|_{H^1}\big)\cdot
\Big[J^\psi(u,v)+J^{\psi^{-1}}
(v,u)\Big]~\leq~ C \cdot J(u,v)\,,\cr}\eqno(3.10)$$
completing the proof of (3.8).
\endproof
\v
\n{\bf Lemma 4.} {\it Let $(u_n)_{n\geq 1}$ be a
Cauchy sequence
for the distance $J$, uniformly
bounded in the $H^1_\per$ norm.
Then
\v
\i{(i)} There exists a limit function
$u\in H^1_\per$ such that $u_n\to u$ in $\L^\infty$
and the sequence of derivatives
$u_{n,x}$ converges to $u_x$ in $L^p_\per$,
for $1\leq p <2$.
\v
\i{(ii)} Let $\mu_n$ be the absolutely continuous measure
having density $u^2_{n,x}$ w.r.t.~Lebesgue measure.
Then one has the weak convergence $\mu_n\wto \mu$, for
some measure $\mu$ whose absolutely continuous part
has density $u_x^2$.
\v
\i{}
}

\n{\bf Proof.}~ {\bf 1.}
By Lemma 3 we already know the convergence
$u_n\to u$, for some limit function $u\in
\L^1_\per\,$.  By a Sobolev embedding theorem,
all functions $u_n,u$  are uniformly H\"older continuous.
This implies $\|u_n-u\|_{\L^\infty}\to 0$.

To establish the convergence of derivatives, we first show that
the sequence of functions
$$v_n\doteq \exp\{ 2i\arctan u_{n,x}\}$$
is compact in $\L^1_\per\,$.
Indeed, fix $\ve>0$.  Then there exists $N$ such that
$J(u_m,u_n)<\ve$ for $m,n\geq N$.
We can now approximate $u_N$ in $H^1_\per$ with a
piecewise affine function
$\tilde u_N$ such that $J(\tilde u_N, u_N)\leq \ve$.
By assumption, choosing suitable transport maps $\psi_n$
we obtain
$$\int_0^1 \Big|\exp \big\{ 2i\arctan u_{n,x}(x)\big\}-
\exp\big\{ 2i\arctan \tilde u_{N,x}(\psi_n(x))\big\}\Big|\,dx
~\leq~ 2\,J(u_n,\tilde u_N)~\leq ~ 4 \ve$$
for all $n\geq N$.
We now observe that all functions
$x\mapsto \exp\big\{ 2i\arctan \tilde u_{N,x}(\psi_n(x))\big\}$
are uniformly bounded, piecewise constant with the same
number of jumps: namely, the number of subintervals
on which $\tilde u_N$ is affine.  The set of all
such functions is compact in $\L^1_\per$.
This argument shows that
the sequence $v_n\doteq \exp\{ 2i\arctan u_{n,x}\}$
eventually remains in an $\ve$-neighborhood of a compact
subset of $\L^1_\per$.  Since $\ve>0$ can be taken arbitrarily
small, by possibly choosing a subsequence we
obtain the strong convergence
$v_n\to v$ for some  $v\in \L^1_\per\,$.
\v
\n{\bf 2.} From the uniform $H^1$ bounds and the $\L^1$ convergence
of the functions $v_n$, we now derive the $\L^p$ convergence of the
derivatives.  For a given $\ve>0$, define
$$M\doteq \sup_n\|u_n\|_{H^1_\per}\,,
\qquad\qquad
A_n\doteq \Big\{ x\in [0,1]\,;~~ \big|u_{n,x}(x)\big| > M/\ve\Big\}\,.$$
The above definitions imply
$$\meas (A_n)\leq \ve^2\eqno(3.11)$$
We now have
$$\eqalign{\|u_{m,x}-u_{n,x}\|_{\L^p}&\leq \left(\int_{A_n\cup A_m}
|u_{m,x}-u_{n,x}|^p\,dx\right)^{1/p}+ \left(\int_{[0,1]\setminus(A_n\cup
A_m)} |u_{m,x}-u_{n,x}|^p\,dx\right)^{1/p}\cr
&\doteq I_1+I_2\,.\cr}\eqno(3.12)$$

$$\eqalign{I_1 &\leq\left[\int_{A_m\cup A_n}
1\cdot dx\right]^{(2-p)/2p} \cdot \left[\int_{A_m\cup
A_n}\big(
|u_{m,x}|+|u_{n,x}|\big)^2\,dx\right]^{1/2}\cr
&\leq \ve^{(2-p)/p}\cdot 2M
\,.\cr}\eqno(3.13)$$

Next, choosing a constant $C_\ve$ such that
$$\big|e^{2i\arctan a}-e^{2i\arctan b}\big|\geq C_\ve
|a-b|\qquad\qquad\hbox{whenever}~~|a|, |b|\leq
M/\ve\,,$$  we obtain
$$I_2 \leq C_\ve \,\left[\int  \Big|e^{2i\arctan u_{m,x}}
-e^{2i\arctan u_{n,x}}\big|^p\,dx\right]^{1/p}.\eqno(3.14)$$
Taking $\ve>0$ small,
we can make the right hand side of (3.13) as small as we like.
On the other hand, choosing a subsequence such that
$v_\nu=e^{2i\arctan u_{\nu,x}}$ converges in $\L^1_\per$,
the right hand side of (3.154) approaches zero.
Hence, for this subsequence,
$$\limsup_{m,n\to\infty}\|u_{m,x}-u_{n,x}\|_{\L^p_\per}=0\,.$$
Since $u_n\to u$ uniformly,
in this case we must have
$$\|u_{n,x}-u_x\|_{\L^p_\per}\to 0\,.\eqno(3.15)$$
We now observe that from any subsequence we can extract a further
subsequence for which (3.15) holds.
Therefore, the whole sequence $(u_{n,x})_{n\geq 1}$
converges to $u_x$ in $\L^p_\per\,$.
\v
\n{\bf 3.}   To establish (ii), we consider
the sequence of measures having density $1+u_{n,x}^2$
w.r.t.~Lebesgue measure.
This sequence converges weakly, because our distance functional is
stronger than the Kantorovich-Waserstein metric which induces
the topology of weak convergence on spaces of measures.
Therefore, $\mu_n\wto \mu$ for some positive measure $\mu$.

Since the sequence $1+u_{n,x}$ converges to $1+u_x$
in $\L^1_\per$, by possibly choosing a subsequence we
achieve the pointwise convergence $u_{n,x}(x)\to u_x(x)$,
for a.e.~$x\in [0,1]$.  For any $\ve>0$, by Egorov's theorem
we have the uniform convergence $u_{n,x}(x)\to u_x(x)$
for all $x\in [0,1]\setminus V_\ve$, for some set
with $\meas(V_\ve)<\ve$.  Since $\ve>0$ can
be taken arbitrarily small, this shows that the absolutely
continuous part of the measure $\mu$ has density $u^2+u_x^2$
w.r.t.~Lebesgue measure.
\endproof

\vsk
\n{\medbf 4 - Continuity in time of the distance functional}
\v
Here and in the next section we examine how the
distance functional $J(\cdot,\cdot)$ evolves in time,
in connection with multi-peakon solutions of
the Camassa-Holm equation (1.1).
We first provide estimates valid on a time interval where no
peakon interactions occur.  Then we show that the
distance functional is continuous across times of interaction.
Since the number of peakons is locally finite, this will suffice
to derive the basic estimates (1.17)-(1.18), in the case of
multi-peakon solutions.
\v
\n{\bf Lemma 5.} {\it Let $t\mapsto u(t)\in H^1_\per$ be a multi-peakon
solution of (1.1).   Assume that no peakon interactions occur
within the interval $[0,\tau]$.  Then
$$J\big(u(s),\, u(s')\big)\leq   C\cdot |s-s'|\,,\qquad\qquad
\forall s,s'\in [0,\tau]\,,\eqno(4.1)$$
for some constant $C$, uniformly valid as $u$ ranges on bounded subsets
of $H^1_\per\,$.}
\v
\n {\bf Proof.}  Assume $0\leq s<s'\leq \tau$.
By the assumptions, the solution $u=u(t,x)$
remains uniformly Lipschitz continuous on the time interval $[0,\tau]$.
Therefore, for each $s\in [0,\tau]$ and $x\in\R$, the Cauchy problem
$${d\over dt}\,\xi(t)= u\big( t,\,\xi(t)\big)\,,\qquad
\qquad \xi(s)=x\,,\eqno(4.2)$$
determines a unique
characteristic curve $t\mapsto \xi(t; s,x)$ passing through the point $(s, x)$.
Given $s'\in [0,\tau]$, we can thus define a transportation plan by setting
$$\psi(x)\doteq \xi(s';s,x)\,.\eqno(4.3)$$
Of course, moving mass along the characteristics
is the most natural thing to do.
We then choose $\phi_1,\phi_2$ to be as large as possible,
according to (3.4).  Namely:
$$\phi_1(x)\doteq
\sup\,\bigg\{\theta\in [0,1]\,;~~ \theta\cdot
\Big( 1+u_x^2(s,x)\Big)\leq \psi'(x)\cdot\Big(
1+ u_x^2
\big(s',\,\psi(x)\big)\Big)\bigg\}\,,$$
$$\phi_2(x)\doteq
\sup\,\bigg\{\theta\in [0,1]\,;~~ \,\theta\cdot
\Big(
1+u_x^2
\big(s',\,\psi(x)\big)\Big)\,\psi'(x)\leq
1+u_x^2(s,x)\bigg\}\,.$$
The cost of this plan is
bounded by
$$\eqalign{J^\psi&\big( u(s)\,,~u(s')\big)\leq
\int_0^1\Big\{ \big|x-\xi(s';s,x)\big|+ \big|
u(s,x)-u(s',\xi(s';s,x))\big|\cr &\qquad\qquad+
\big|2\arctan u_x(s,x)-2\arctan  u_x(s',\xi(s';s,x))
\big|_*\Big\}\big(1+u_x^2(s,x)\big)\,dx\cr
&\qquad +\int_0^1 \big(1-\phi_1(x)\big)\,\big(1+u^2_x(s,x)\big)\,dx\cr
&\qquad
+\int_0^1 \big(1-\phi_2(\psi(x))\big)\,\Big(1
+u^2_x\big(s',\psi(x)\big)\Big)\,\psi'(x)\,dx
\,.
\cr}\eqno(4.4)$$

To estimate the right hand side of (4.4), we first
observe that, for all $u\in H^1_\per$,
$$\eqalign{\big\|u\big\|_{\L^\infty}
&\leq \int_0^1 \big|u(x)\big|\,dx+\int_0^1
\big|u_x(x)\big|\,dx\cr
&\leq \|u\|_{\L^2}+\|u_x\|_{\L^2}\leq 2\|u\|_{H^1_\per}=2(E^u)^{1/2}
.\cr}\eqno(4.5)$$
Using (4.5) in (4.2) we obtain
$$\big|\xi(s)-\xi(s')\big|\leq 2(E^{\bar u})^{1/2}\cdot |s-s'|\,.\eqno(4.6)$$

Next, from the definition of the source term $P$ at (1.2)
it follows
$$\big\|P\|_{\L^\infty}\leq {1\over 2} \big\|e^{-|x|}\,\big\|_{\L^1(\R)}
\cdot \left\| u^2+{u_x^2\over 2}\right\|_{\L^1([0,1])}\leq
\|u\|^2_{H^1_\per}
=E^u\,.\eqno(4.7)$$
Similarly,
$$\big\|P_x\|_{\L^\infty}\leq \|u\|^2_{H^1_\per}
=E^u\,.\eqno(4.8)$$
Using (4.8) we obtain
$$\eqalign{\Big| u\big(s',\,\xi(s')\big)&-u\big(s,\,\xi(s)\big)\Big|
\leq \int_s^{s'}\left| {d\over dt} u\big(t,\,\xi(t)\big)\right|\,dt\cr
&=\int_s^{s'}\Big| P_x\big(t,\,\xi(t)\big)\Big|\,dt~\leq
~E^{\bar u} \cdot |s'-s|\,. \cr}\eqno(4.9)$$

Concerning the term involving arctangents,
recalling (1.3) we obtain
$${d\over dt}\Big[2\arctan u_x\big(t,\xi(t,x)\big)\Big]={2\over 1+u_x^2}\,
\left[ u^2-{u_x^2\over 2}-P\right]\,.$$
The bounds (4.7) and (4.9) thus yield
$$\eqalign{\Big|2\arctan u_x\big(s',\xi(s')\big)
-2\arctan u_x\big(s,\,\xi(s)\big)\Big|_*
&\leq \Big( 2\|u\|^2_{\L^\infty} +1+2\|P\|_{\L^\infty}\Big)\cdot |s'-s|
\cr
&\leq \big(10\,E^{\bar u}+1\big)\cdot |s'-s|\,.
\cr}
\eqno(4.10)$$
This already provides a bound on
the first integral on the right hand side of (4.4).

Next, call $I_1,I_2$ the last two integrals on the right hand side of (4.4).
Notice that
$$I_1+I_2=\int_0^1 \bigg|\Big(1+ u_x^2(s,y)\Big)-\xi_y(s';s,y)\,
\Big( 1+ u_x^2\big(s',\xi(s';s,y)\big)\Big)\bigg|\,dy\eqno(4.11)$$
Indeed, $I_1+I_2$ measures the difference between
the measure $\big( 1+u_x^2(s',y)\big)\,dy$ and the push-forward
of the measure $\big( 1+u_x^2(s,x))\,dx$ through the mapping
$x\mapsto \xi(s';s,x)$.

Since the push-forward of the measure $u_x^2\,dy$ satisfies the
linear conservation law
$$w_t+(uw)_x=0\,,\eqno(4.12)$$
comparing (4.12) with (1.5) we deduce
$$\eqalign{
\int_0^1 \Big| &u_x^2(s,y)-\xi_y(s,y)\,
 u_x^2\big(s',\xi(s';s,y)\big)\Big|\,dy~\leq
~\int_s^{s'} \int_0^1 2\big| (u^2-P)u_x\big|\,dx\,dt\cr
&\leq
2 \int_s^{s'}  \big( \|u\|_{\L^\infty}^2+
\|P\|_{\L^\infty}\big)\,\|u_x\|_{\L^1}\,dt
~\leq~ 2\,\big(4E^{\bar u}+E^{\bar u}\big)\,E^{\bar u}\cdot |s'-s|\,,
\cr
}
\eqno(4.13)$$
because of (4.5), (4.7) and (1.15).
Finally, we need to estimate the remaining terms,
describing by how much the Lebesgue measure fails to be conserved
by the transformation  $x\mapsto \xi(s';s,x)$.
Observing that
$${\partial \over\partial t}\,\xi_y(t,y)= u_x\big(t,\,\xi(t,y)\big)\,
\xi_y(t,y)\,,\qquad\qquad\xi_y(0,y)= 1\,,\eqno(4.14)$$
we find
$$\eqalign{\int_0^1
\big|1&  -\xi_y(s';s,y)\,
\big|\,dy\leq
\int_s^{s'} \int_0^1\bigg|{\partial \over \partial t}\big[
\xi_y(t,s,y)\big]\bigg|
\,dy\,dt\cr
&\leq \int_s^{s'}\int_0^1
\xi_y(t;s,y)\,
\big|u_x(t,\xi(t;s,y))\big|
\,dy\,dt
\,.\cr}\eqno(4.15)$$
To estimate the right hand side of (4.15),
we use the decomposition
$[0,1]=Y\cup Y'\cup Y''$, where
$$Y\doteq \Big\{ y\,;~~ \xi_y(t;s,y)\in \big[(1/2)\,,~2\big]\quad
\hbox{for all}~t\in [s,s']\,\Big\}\,,$$
$$Y'\doteq \Big\{ y\,;~~ \xi_y(t;s,y) <1/2 \quad
\hbox{for some}~t\in [s,s']\,\Big\}\,,$$
$$Y''\doteq \Big\{ y\,;~~ \xi_y(t;s,y) > 2\quad
\hbox{for some}~t\in [s,s']\,\Big\}\,.$$
Integrating over $Y$ one finds
$$
\int_s^{s'} \int_Y
\xi_y(t;s,y)
\,\big|u_x(t,\xi(t;s,y))\big|
\,dy\,dt~\leq~ 2\int_s^{s'}
\big\|u_x(t)\big\|_{\L^1}\,dt~\leq~ 2 E^{\bar u}\cdot |s'-s|\,.
\eqno(4.15)$$
Next, if $y\in Y'$ we define
$$\tau(y)\doteq \inf\,\big\{ t>s\,;~~\xi_y(t;s,y)< 1/2\big\}\,.$$
Observe that $y\in Y'$ implies
$$\int_s^{\tau(y)} \Big|u_x\big(t,\,\xi(t;s,y)\big)\Big|\,dt
\geq \ln 2\,.$$
Therefore
$$\eqalign{\int_{Y'} dy&\leq
{1\over \ln 2}\int_{Y'}
\left[\int_s^{\tau(y)}
\Big|u_x\big(t,\,\xi(t;s,y)\big)\Big|\,dt
\right]\,dy\cr
&\leq
{2\over \ln 2}\int_s^{s'}\int_0^1
\big|u_x\big(t,x)\big|\,dx\,dt~\leq~ 4 \,E^{\bar u} \cdot |s'-s|\,.
\cr}\eqno(4.16)$$
The estimate for the integral over $Y''$ is entirely analogous,
Indeed, the push-forward of the Lebesgue measure along characteristic
curves from $t=s$ to
$t=s'$ satisfies exactly the same type of estimates as
the pull-back of the Lebesgue measure from $t=s'$ to
$t=s$. All together, these three estimates imply
$$\int_{Y\cup Y'\cup Y''}
\big|1 -\xi_y(s';s,y)\,
\big|\,dy\leq 10\,E^{\bar u}\cdot |s'-s|\eqno(4.17)$$
\v
Putting together the estimates (4.6), (4.9), (4.10), (4.3) and (4.17),
the distance in (4.4) can be estimated by
$$
J^\psi\big( u(s)\,,~u(s')\big)\leq\Big[
2(1+E^{\bar u}) +E^{\bar u} + (10 \,E^{\bar u} +1)+
10\,(E^{\bar u})^2+10\,E^{\bar u}\Big]\cdot |s'-s|\,.
\eqno(4.18)$$
This establishes (4.1).
\endproof
\v
According to Lemma 5, as long as no peakon interactions occur,
the map $t\mapsto u(t)$ remains
uniformly Lipschitz continuous w.r.t.~our distance functional,
with a Lipschitz constant that depends only
on the total energy $E^{\bar u}$.  Since interactions
can occur only at isolated times, to obtain a global Lipschitz estimate
it suffices to show that trajectories are continuous
(w.r.t.~the distance $J$) also at interaction times.
\v
\n{\bf Lemma 6.} {\it Assume that the multi-peakon solution
$u(\cdot)$ contains two or more peakons which interact at
a time $\tau$.   Then
$$\lim_{h\to 0+} J\big( u(\tau -h),\, u(\tau + h)\big)=0\,.\eqno(4.19)$$
}\v
\n{\bf Proof.} To fix the ideas, call $x=\bar q$ the place
where the interaction occurs,
and let $p_1,\ldots, p_k$ be the strengths of the peakons
that interact at time
$\tau$.  We here assume that $0<\bar q<1$.
The case where two or more groups of peakons interact
exactly at the same time $\tau$, within the interval
$[0,1]$, can be treated similarly.

For $|t-\tau|\leq h$, call $\xi^-(t)$, $\xi^+(t)$ respectively
the position of the smallest and largest characteristic
curves passing through the point $(\tau,\bar q)$, as in (2.5).
We observe that $u$ is Lipschitz continuous
in a neighborhood of each point $(\tau,x)$, with $x\not= \bar q$.
Hence, for $x\in [0,1]\setminus\{\bar q\}$ there exists a unique
characteristic curve $t\mapsto \xi(t; \tau,x)$ passing through
$x$ at time $\tau$.
For a fixed $h>0$, the transport map $\psi$ is defined as follows.
Consider the intervals
$I_{-h}\doteq \big[\xi^-(\tau-h),\,\xi^+(\tau-h)\big]$
and $I_h\doteq \big[\xi^-(\tau+h),\,\xi^+(\tau+h)\big]$.
On the complement $[0,1]\setminus I_{-h}$
we define
$$\psi\big(\xi(\tau-h;\tau ,x)\big)=\xi(\tau+h;\tau ,x)\,,
$$
so that transport is performed along characteristic curves.
It now remains to extends $\psi$ as a map from $I_{-h}$ onto $I_h$.
Toward this goal, we recall that
our construction of multi-peakon solutions in Section 2
was specifically designed in order to achieve the identity
$$e_{(\tau,\bar q)}=\lim_{h\to 0+} \int_{I_{-h}} u_x^2(\tau-h,x)\,dx~=~
\lim_{h\to 0+} \int_{I_h} u_x^2(\tau+h,x)\,dx\,.\eqno(4.20)$$
For $h>0$ we introduce the quantities
$$E(-h)\doteq \int_{I_{-h}} \big(1+u_x^2(\tau-h,\,x)\big)\,dx\,,\qquad\qquad
E(h)\doteq \int_{I_h} \big(1+u_x^2(\tau+h,\,x)\big)\,dx\,,$$
$$e(h)\doteq 2\min\big\{E(-h),\,E(h)\big\}-\max\big\{E(-h),\,E(h)\big\}\,.$$
Notice that (4.20) implies $e(h)>0$ and
$$ E(-h)\to e_{(\tau,\bar q)}\,,
\qquad E(h)\to e_{(\tau,\bar q)}\,,
\qquad e(h)\to e_{(\tau,\bar q)}\,,\eqno(4.21)$$
as $h\to 0+$.
Consider the point $x^*=x^*(h)$ inside the
interval $I_{-h}=\big[\xi^-(\tau-h),\,\xi^+(\tau-h)\big]$,
implicitly defined by
$$\int_{\xi^-(\tau-h)}^{x^*}\big(1+u_x^2(x)\big)\,dx = e(h)\,.$$
For $x\in \big[\xi^-(\tau-h),~x^*\big]$ we define $\psi(x)$ as the
unique point such that
$$\int_{\xi^-(\tau+h)}^{\psi(x)}\big(1+u_x^2(\tau+h,\,x)\big)\,dx
=\int_{\xi^-(\tau-h)}^x\big(1+u_x^2(\tau-h,\,x)\big)\,dx\,.\eqno(4.22)$$
We then extend $\psi$ as an affine map from
$\big[x^*,\,\xi^+(\tau-h)\big]$ onto $\big[\psi(x^*),\,\xi^+(\tau+h)\big]$,
namely
$$\psi\Big(\theta\cdot\xi^+(\tau-h)+(1-\theta)
\cdot x^*\Big)=\theta\cdot\xi^+(\tau+h)+(1-\theta)
\cdot \psi(x^*)\qquad\qquad \theta\in [0,1]\,.$$
Finally, we prolong $\psi$ to the whole real line according to
(3.3).
As usual, the 1-periodic
functions $\phi_1,\phi_2$ are then chosen to be as large as possible,
according to (3.4).
As $h\to 0+$, we claim that the following quantity approaches zero:
$$\eqalign{&J^\psi\big(u(\tau-h),~u(\tau+h)\big)\cr
&=
\int_0^1 d^\diamondsuit\Big( \big(x,\,u(\tau-h,x),\,
2 \arctan u_x(\tau-h,x)\big)\,,~
\big(\psi(x),\,u(\tau+h, \,\psi(x)),\,2 \arctan \tilde u_x(\tau+h,\,
\psi(x))\big)\Big) \cr
&\qquad\qquad\qquad \cdot
\phi_1(x)\,\big(1+u_x^2(\tau-h,\,x)\big)\,dx\cr
&\qquad+\int_0^1 \big(1-\phi_1(x)\big)\,
\big( 1+u_x^2(\tau-h,\,x)\big)\,dx\cr
& \qquad +\int_0^1
\big(1-\phi_2(\psi(x))\big)\,\Big(1+u_x^2
\big(\tau+h,\,\psi(x)\big)\Big)\,\psi'(x)\,dx
\,.\cr}\eqno(4.23)$$
It is clear that the restriction of all the above
integrals to the complement
$[0,1]\setminus I_{-h}$ approaches zero as $h\to 0$.
We now prove that their restriction to $I_{-h}$ also vanishes
in the limit.
As $h\to 0+$, for $x\in I_{-h}$ we have
$$d^\diamondsuit\Big(
\big(x,\,u(\tau-h,x),\,2\arctan u_x(\tau-h,x)\big)~,~
\big(\psi(x),\,u(\tau+h,\psi(x)),\,2\arctan u_x(\tau+h,\psi(x))\big)\Big)
~\to~ 0\,,$$
because all points approach the same limit
$\big(\bar q, \,u(\tau,\bar q),\, \pi\big)$.
The first integral in (4.23) thus approaches zero as $h\to 0+$.

Concerning the last two integrals, by (4.22) it follows
$$\phi_1(x)=\phi_2\big(\psi(x)\big)=1\qquad\qquad \forall x\in
\big[\xi^-(\tau-h),\,x^*\big]\,.$$
Moreover, our choice of $x^*$ implies
$$\eqalign{
&\int_{x^*}^{\xi^+(\tau-h)}
\big( 1+u_x^2(\tau-h,\,x)\big)\,dx+\int_{\psi(x^*)}^{\xi^+(\tau+h)}
\Big(1+u_x^2
\big(\tau+h,\,\psi(x)\big)\Big)\,\psi'(x)\,dx\cr
&\qquad \leq
2\max\big\{E(-h),\,E(h)\big\}-2\min\big\{E(-h),\,E(h)\big\}\,.
\cr}\eqno(4.24)$$
By (4.21), as $h\to 0+$ the
right hand side of (4.24) approaches zero.
Hence the same holds for the last two integrals in (4.23).
This completes the proof of the lemma.
\endproof
\vsk

\n{\medbf 5 -
Continuity w.r.t.~the initial data}
\v
We now consider two distinct
solutions and study how the distance $J\big(u(t)\,~
v(t)\big)$ evolves in time.
To fix the ideas, let $t\mapsto u(t)$ and $t\mapsto v(t)$
be two multi-peakon solutions of (1.1), and assume that
no interaction occurs within a given time interval $[0,T]$.
In this case, the functions $u,v$ remain Lipschitz continuous.
We can thus define the characteristic curves
$t\mapsto \xi(t,y)$ and $t\mapsto \zeta(t,\tilde y)$
as the solutions to the Cauchy problems
$$\eqalign{\dot \xi &= u(t,\xi),\qquad\qquad \xi(0)=y\,,\cr
\dot \zeta &= v(t,\zeta),\qquad\qquad \zeta(0)=\tilde y\,,\cr}\eqno(5.1)$$
respectively.
Let now $\psi_{(0)}\in\F$ be any transportation plan at time $t=0$.
For each $t\in [0,T]$ we can define a transportation plan
$\psi_{(t)}\in\F$ by setting
$$\psi_{(t)}\big(\xi(t,y)\big)\doteq \zeta\big(t,~\psi_{(0)}
(y)\big).\eqno(5.2)$$
The corresponding functions
$\phi_{1}^{(t)},\phi_2^{(t)}$
are then defined according to (3.4), namely
$$\phi_1^{(t)}(x)\doteq
\sup\,\bigg\{\theta\in [0,1]\,;~~ \theta\cdot
\Big( 1+u_x^2(t,x)\Big)~\leq ~\Big(
1+ v_x^2
\big(t,\,\psi_{(t)}(x)\big)\Big)\,\psi_{(t)}'(x)\bigg\}\,,$$
$$\phi_2^{(t)}(x)\doteq
\sup\,\bigg\{\theta\in [0,1]\,;~~ 1+u_x^2(t,x)~\geq ~\theta\cdot
\Big(
1+v_x^2
\big(t,\,\psi_{(t)}(x)\big)\Big)\,\psi'_{(t)}(x)\bigg\}\,.$$

If
initially the point $y$ is mapped into $\tilde y\doteq \psi_{(0)}(y)$,
then at a later time $t>0$ the point $\xi(t,y)$ along the
$u$-characteristic starting from $y$ is sent to the point  $\zeta
(t,\tilde y)$ along the $v$-characteristic
starting from $\tilde y$. We thus transport
mass from the
point $\Big(\xi(t,y)\,,~u\big(t,\xi(t,y)\big)\,,~ 2\arctan
u_x\big(t,\xi(t,y)\big)\Big)$ to the corresponding point
$\Big(\zeta(t,\tilde y)\,, ~v\big(t,\zeta
(t,\tilde y)\big)\,, ~2\arctan v_x\big(t,\zeta
(t,\tilde y)\big)\Big)\,,$

In the following, our main goal is to provide an upper bound
on the
time derivative of the function
$$\eqalign{
J^{\psi(t)}&\big(u(t)\,,~ v(t)\big)
\doteq
\int_0^1 d^\diamondsuit\Big( \big(x,\,u(t,x),\,
2 \arctan u_x(t,x)\big)\,,\cr
&\qquad\qquad \big(\psi_{(t)}(x),\,v(t,\,\psi_{(t)}(x)),\,
2 \arctan v_x(t,\,\psi_{(t)}(x))\big)\Big) \cdot
\phi_1^{(t)}(x)\,\big(1+u_x^2(t,x)\big)\,dx\cr
&+\int_0^1 \big(1-\phi_1^{(t)}(x)\big)\,
\big( 1+u_x^2(t,x)\big)\,dx\cr
& +\int_0^1
\big(1-\phi_2^{(t)}(\psi_{(t)}(x))\big)\,\Big(1+v_x^2
\big(t,\,\psi_{(t)}(x)\big)\Big)\,\psi_{(t)}'(x)\,dx
\,.\cr}\eqno(5.3)$$
Differentiating the right hand side of (5.3)
one obtains several terms, due to
\v
\item{$\bullet$}~Changes in the distance $d^\diamondsuit$
between the points
$(\xi,\, u,\, 2\arctan u_x)$ and $(\zeta,\, v, \,
2\arctan v_x)$.
\v
\item{$\bullet$}~Changes in the base measures $(1+u_x^2)\,dx$ and
$(1+v_x^2)\,dx$.
\v
\n Throughout the following, by $\O(1)$ we denote a quantity which remains
uniformly bounded as $u,v$ range in bounded subsets of $H^1_\per\,$.
Using the elementary estimate
$$|u-v|\leq \big(1+|u|+|v|\big)\,\min\big\{ |u-v|,\,1\big\}\,,$$
we
begin by deriving the bound
$$\eqalign{I_1&\doteq \int_0^1 {d\over dt}\,\big|  x-\psi_{(t)}(x)
\big|\cdot\phi_1^{(t)}(x)\,\big( 1+u_x^2(t,x)\big) \,\,dx\cr
&\leq \int_0^1 \Big| u(t,x)-v\big(t,\,\psi_{(t)}(x)\big)
\Big|\cdot \phi_1^{(t)}(x)\,\big( 1+u_x^2(t,x)\big) \,\,dx\cr
&\leq \Big(1+ \big\|u(t)\big\|_{\L^\infty} +
\big\|v(t)\big\|_{\L^\infty}
\Big)\cdot J^{\psi_{(t)}}\big(u(t),\,v(t)\big)\cr
&= \O(1)\cdot J^{\psi_{(t)}}\big(u(t),\,v(t)\big)\,.
}\eqno(5.4)
$$
Here and in the sequel, the time derivative
is computed along characteristics.

Next,
recalling the basic equation (1.1), we consider
$$\eqalign{I_2&\doteq \int_0^1 {d\over dt}\Big| u(t,x)- v(t,\,
\psi_{(t)}(x)
\big)\Big|\cdot\phi_1^{(t)}(x)\,\big( 1+u_x^2(t,x)\big) \,\,dx\cr
&\leq\int_0^1 \Big| P_x^u(t,x)-P_x^v\big(t,\psi_{(t)}(x)\big)
\Big|\cdot\big( 1+u_x^2(t,x)\big) \,\,dx\,.\cr}
\eqno(5.5)$$
In the spatially periodic case, by (1.2) and (2.2) we
can write the source terms $P^u_x$, $P^v_x$ as
$$\eqalign{ P^u_x(t,x)&={1\over 2}\int_{x-1}^x \chi'(x-y)\cdot
\left[u^2(t,y)+{u_x^2(t,y)\over 2}\right]dy\,,
\cr
P^v_x\big(t,\,\psi_{(t)}(x)\big)&=
{1\over 2}\int_0^1 \chi'\big(
\psi_{(t)}(x)-\tilde y\big)
\cdot v^2(t,\tilde y)\,d\tilde y\,\cr
&\qquad\qquad +\int_{x-1}^x \chi'\big(
\psi_{(t)}(x)-\psi_{(t)}(y)\big)
\cdot {v_x^2\big(t,\psi_{(t)}(y)\big)
\over 4}\psi'_{(t)}(y)\,dy\,,\cr
}\eqno(5.6)$$
where, according to (2.2),
$$\chi'(x)= { e^x-e^{1-x}\over e-1}\qquad 0<x<1\,,\qquad\qquad
\chi'(x)=\chi'(x+1)\qquad x\in\R\,.\eqno(5.7)$$

%In the following we assume that the transportation plan satisfies
%$$\psi'_{(t)}(x)=\tilde \psi'_{(t)}(x)
%\doteq{1+u_x^2(t,x)\over 1+v_x^2(t,\psi_{(t)}(x))}\,.\eqno(5.7)$$

In the next computation, we use the estimate
$$\int_0^1 \bigg|\big( 1+u_x^2(y)\big)- \Big( 1+v_x^2\big(\psi(y)\big)
\Big)\psi'(y)\bigg|\,dy=\O(1)\cdot J^\psi (u,v)\,.\eqno(5.8)$$
which holds because of the last two terms in the definition (3.6).
Observing that $\chi'$ is Lipschitz continuous on the open interval
$]0,1[\,$, we now compute (omitting explicit references to the time
$t$)
$$
\eqalign{ \Big|P^u_x(x) &-P^v_x\big(\psi(x)\big)\Big|
\leq
{1\over 2} \int_0^1 \left|\chi'(x-y)\cdot u^2(y)-\chi'\big(\psi(x)-y\big)
\cdot v^2(y)\right|\,dy
\cr
&+\O(1)\cdot\int_{x-1}^x \Big|x-y-(\psi(x)-\psi(y)\big)\Big|\cdot
{v_x^2\big(\psi(y)\big)\over 2}\,
\psi'(y)\,dy \cr
&
+{1\over 4}\left|
\int_{x-1}^x
\chi'(x-y)\,\Big(u_x^2(y)-v_x^2\big(\psi(y)\big)\,\psi'(y)\Big)\,dy
\right|\cr
&= \O(1)\cdot \Big(\big|x-\psi(x)\big|+ \|u^2-v^2\|_{\L^1}\Big)\cr
&\qquad
+\O(1)\cdot \left(\big|x-\psi(x)\big|+ \int_{x-1}^x \big|y-\psi(y)\big|\cdot
{v_x^2\big(\psi(y)\big)\over 2}\,
\psi'(y)\,dy \right)
\cr
&\qquad
+\O(1)\cdot \left(J^\psi (u,v)+
\Big|\int_{x-1}^x \chi'(x-y)\cdot \big[\psi'(y)-1\big]
\,dy\Big|\right)
\cr
&= \O(1)\cdot \big|x-\psi(x)\big|+ \O(1)\cdot J^\psi (u,v)+
\O(1)\cdot \left( \big|x-\psi(x)\big|+\int_{x-1}^x
\chi''(x-y)\cdot \big[\psi(y)-y\big]
\,dy\right)
\cr
&=\O(1)\cdot \big|x-\psi(x)\big|+ \O(1)\cdot J^\psi (u,v)\,.\cr}
\eqno(5.9)
$$
Integrating over one period we conclude
$$I_2 = \O(1)\cdot J^\psi \big(u(t), \,v(t)\big)\,.\eqno(5.10)$$

For future use, we observe that a computation entirely similar to
(5.9) yields
$$\Big|P^u_x(x) -P^v_x\big(\psi(x)\big)\Big|
=\O(1)\cdot \big|x-\psi(x)\big|+ \O(1)\cdot J^\psi (u,v)\,.\eqno(5.11)$$

Next, we look at the term
$$I_3\doteq \int_0^1 {d\over dt}
\Big| 2\arctan\, u_x(t,x)
-2\arctan\, v_x\big(t,\psi_{(t)} (x)\big)
\Big|\cdot\phi_1^{(t)}(x)\,\big( 1+u_x^2(t,x)\big)\,dx\,.\eqno(5.12)$$
Along a characteristic, according to (1.3) one has
$${d\over dt} 2\arctan u_x(t,\xi(t))={2\over 1+u_x^2}\Big[u^2-
{u_x^2\over 2}-P^u\Big]\,.$$
Call
$\theta^u\doteq 2\arctan u_x\,$,
$\theta^v\doteq 2\arctan v_x\,$,
so that
$${1\over 1+u_x^2}=\cos^2{\theta^u\over 2}\,,\qquad\quad
{u_x\over 1+u_x^2}={1\over 2}\sin \theta^u\,,\qquad\quad {u_x^2
\over 1+u_x^2}=\sin^2{\theta^u\over 2}\,.$$
We now have
$$\eqalign{\int_0^1 &\big(1+u_x^2(x)\big)\cdot
\left| {u_x^2(x)\over 1+u_x^2(x)}- {v_x^2(\psi(x))
\over 1+v_x^2(\psi(x))}\right|\,dx\cr
&
=
\int_0^1 \big(1+u_x^2(x)\big)\cdot
\Big| \sin^2{\theta^u(x)\over 2}-\sin^2{\theta^v\big(\psi(x)\big)\over 2}
\Big|
\,dx\cr
&\leq \int_0^1 \big(1+u_x^2(x)\big)\cdot \Big|\theta^u(x)-\theta^v
\big(\psi(x)\big)\Big|\,dx\cr
&=\O(1)\cdot J(u,v)\,.\cr}\eqno(5.13)$$
Next, using (5.11) we compute
$$\eqalign{
\int_0^1\big(1+u_x^2(x)\big)\cdot &\bigg|
{u^2(x)-P^u(x)\over 1+u_x^2(x)}-{ v^2\big(\psi(x)\big)-P^v\big(\psi(x)\big)
\over 1+v_x^2\big(\psi(x)\big)}\bigg|\,dx\cr
&\leq
\int_0^1
\Big|u^2(x)-v^2\big(\psi(x)\big)\Big|\,dx
+\int_0^1
\Big|P^u(x)-P^v\big(\psi(x)\big)
\Big|\,dx\cr
&\qquad\qquad+\O(1)\cdot \int_0^1 \left| {1\over 1+u_x^2(x)}-
{1\over 1+v_x^2\big(\psi(x)\big)}\right|\cdot \big(1+u_x^2(x)\big)\,dx
\cr
&=\O(1)\cdot J^\psi(u,v)\,,\cr}
\eqno(5.14)$$
where the last term was estimated by observing that
$$\left|{1\over 1+u_x^2}-
{1\over 1+v_x^2}\right|\leq
\big|2\arctan u_x-2\arctan v_x\big|_*
\,.\eqno(5.15)$$
Putting together all previous estimates we conclude
$$I_1+I_2+I_3=\O(1)\cdot J^\psi (u,v)\,.\eqno(5.16)$$
%This shows that the rate of change in $J^{\psi_{(t))} \big(u(t),\,v(t)
%\big)$ due to changes in the distance $d^\diamondsuit$

\vs

To complete the analysis, we have to consider the terms due to the change
in base measures. From (1.5) it follows that
the production of new mass in the base measures is described
by the balance laws
$$\left\{
\eqalign{(1+u_x^2)_t+ \big[ u(1+u_x^2)\big]_x&= [2u^2+1-2P^u]u_x\doteq
f^u\,,\cr
(1+v_x^2)_t+ \big[ v(1+v_x^2)\big]_x&= [2v^2+1-2P^v]v_x \doteq
f^v\,.\cr}\right.\eqno(5.17)$$
This leads us to consider two further integrals $I_4,I_5$ :
$$\eqalign{ I_4&=\int_0^1
 d^\diamondsuit\Big( \big(x,\,u(x),\, 2\arctan u_x(x)\big)\,,~
\big(\psi(x),\,v(\psi(x)),\,2\arctan v_x(\psi(x))\Big) \cr
&\qquad\qquad \cdot
\big|2u^2(x)+1-2P^u(x)\big|\,\big|u_x(x)\big|\,dx\cr
&=\O(1)\cdot
\int_0^1
d^\diamondsuit\Big( \big(x,\,u(x),\, 2\arctan u_x(x)\big)\,,~
\big(\psi(x),\,v(\psi(x)),\,2\arctan v_x(\psi(x))\Big) \cr
&\qquad\qquad \qquad\qquad\cdot
\big(1+u_x^2(x)\big)\,dx\cr
&=\O(1)\cdot J^\psi(u,v)\,.\cr}\eqno(5.18)$$

$$\eqalign{ I_5&=\int_0^1\bigg|
\big[2u^2(x)+1-2P^u(x)\big]u_x(x)-\Big[ 2v^2\big(\psi(x)\big)
+1-2P^v\big(\psi(x)\big)
\Big] v_x\big(\psi(x)\big)\,\psi'(x)\bigg|\,dx
\cr
&\leq 2
\int_0^1\Bigg\{\Big| u^2(x)-v^2\big(\psi(x)\big)\Big|+\Big|P^u(x)-
P^v\big(\psi(x)\big)
\Big|\Bigg\}\, \big|u_x(x)\big|\,dx
\cr
&\qquad
+\int_0^1 \Big|2v^2\big(\psi(x)\big)
+1-2P^v\big(\psi(x)\big)
\Big|\cdot \Big|u_x(x)-v_x\big(\psi(x)\big)
\psi'(x)
\Big|\,dx
\cr
&=I_5'+I_5''\cr
}\eqno(5.19)
$$

Using (5.11) we easily obtain
$$\eqalign{
I_5'&\leq \int_0^1\Bigg\{\Big| u^2(x)-v^2\big(\psi(x)\big)
\Big|+\Big|P^u(x)-
P^v\big(\psi(x)\big)
\Big|\Bigg\}\, \big(1+u_x^2(x)
\big)\,dx\,\cr
&=\O(1)\cdot J^\psi(u,v)\,.\cr}\eqno(5.20)$$
On the other hand, recalling (5.13) and using the change of variable
$y=\psi(x)$, $x=\psi^{-1}(y)$, we find
$$\eqalign{I_5''&=
\O(1)\cdot
\int_0^1 \Big|u_x(x)-v_x\big(\psi(x)\big)
\psi'(x)
\Big|\,dx\cr
&=\O(1)\cdot\int_0^1  \left|{u_x(x)\over 1+u_x^2(x)}-{
v_x\big(\psi(x)\big)\over 1+v_x^2\big(\psi(x)\big)}
\right|\big(1+u_x^2(x)\big)\,dx\cr
&\qquad
+\O(1)\cdot\int_0^1
\Big|v_x\big(\psi(x)\big)\,\psi'(x)\Big|\cdot\left|{1+u_x^2(x)
\over \big(1+v_x^2(\psi(x))\big)\,\psi'(x)}-1\right|
\cr
&\leq \O(1)\cdot J^\psi(u,v)
+\int_0^1
\Big|\big(1+u_x^2(x)\big)-
\big(1+v_x^2(\psi(x))\big)\,\psi'(x)\Big|\,dx
\cr
&=\O(1)\cdot J(u,v)\,.\cr}\eqno(5.21)$$

All together, the previous estimates show that
$${d\over dt} J^{\psi_{(t)}}\big( u(t),\,v(t)\big)\leq
I_1+I_2+I_3+I_4+I_5'+I_5'' =\O(1)\cdot
J^{\psi_{(t)}}\big( u(t),\,v(t)\big)\,,\eqno(5.22)$$
where $\O(1)$ denotes a quantity which remains uniformly bounded
as $u,v$ range on bounded sets of $H^1_\per\,$.
As an immediate consequence we obtain
\v
\n{\bf Lemma 7.} {\it Let $t\mapsto u(t)$, $t\mapsto v(t)$
be two conservative, spatially periodic
multipeakon solutions, as in Lemma 2.
Then there exists a constant $\kappa$, depending only on
$\max\big\{\|u\|_{H^1_\per}\,,~\|v\|_{H^1_\per}\big\}$, such that
$$J\big( u(t)\,,~v(t)\big)\leq e^{\kappa|t-s|}\cdot
J\big( u(s)\,,~v(s)\big)\qquad\qquad s,t\in\R\,.
\eqno(5.23)$$
}
\v
\n{\bf Proof.}
For $t>s$ the estimate (5.23) follows
from (5.22), taking the infimum among all transportation
plans $\psi_{(s)}$ at time $s$.   The case $t<s$
is obtained simply by observing that the Camassa-Holm
equations are time-reversible.
\endproof

\vsk

\n{\medbf 6 - Proof of the main theorems
}
\v
Thanks to the analysis in the previous sections,
we now all the ingredients
toward a proof of Theorems 1 and 2.

The estimates in (1.16) follow from Lemma 3.
Given an initial data $\bar u\in H^1_\per$,
to construct the solution
of the Camassa-Holm equation
we consider a sequence of multi-peakons $\bar u_n$,
converging to $\bar u$ in $H^1_\per$.
Then we consider the corresponding solutions
$t\mapsto u_n(t)$, defined for all $n\geq 1$ and $t\in\R$.
This is possible because of Lemmas 1 and 2.

We claim that the sequence $u_n(t)$ is Cauchy in $L^2_\per\,$.
Indeed, by Lemma 3 and Lemma 5,
$$\big\|u_m(t)-u_n(t)\big\|_{\L^1_\per}
\leq C\cdot J\big(u_m(t),\,u_n(t)\big)\leq
C\cdot e^{\kappa |t|} \,J\big(u_m(0),\,u_n(0)\big)
\leq C^2\cdot e^{\kappa |t|}\big\|u_m(t)-u_n(t)\big\|_{H^1_\per}\,.$$
Therefore, $u_n(t)\to u(t)$ in $\L^1_\per$, for some function
$u:\R\mapsto H^1_\per\,$.   By interpolation, the convergence
$u_n\to u$ also holds in all spaces $\L^p_\per$, $1\leq p\leq\infty$.
The continuity estimates (1.17)-(1.18) now follow by passing to
the limit in Lemma 5 and 6.

It remains to show that the limit function $u(\cdot)$
is actually a solution to the Camassa-Holm equation and
its energy $E(t)$ in (1.15) is a.e.~constant.
Toward these goals, we observe that all solutions $u_n$
are Lipschitz continuous with the same Lipschitz
constant, as maps from $\R$ into $ L^2_\per\,$. Indeed
$$\|u_{n,t}\|_{\L^2_\per}\leq \|u_n\|_{\L^\infty}\cdot
\|u_{n,x}\|_{\L^2_\per }+\left\| {1\over 2}\,e^{-|x|}\right\|_{\L^2}
\cdot \left\|u_n^2+{u^2_{n,x}\over 2}\right\|_{\L^1_\per}.
\eqno(6.1)$$
As a consequence, the map $t\mapsto u(t)$ has uniformly bounded
$H^1_\per$ norm, and is Lipschitz continuous with values
in $L^2_\per\,$.   In particular,
$u$ is uniformly H\"older continuous as a function of $t,x$
and the convergence $u_m(t,x)\to u(t,x)$ holds uniformly
for $t$ in bounded sets.
Moreover, since $\L^2_\per$ is a reflexive space,
the time derivative $u_t(t)\in\L^2_\per$ is well
defined for a.e.~$t\in\R$.

We now observe that, for each $n\geq 1$,
both sides of the equality
$${d\over dt} u_n=-u_n\, u_{n,x} -P^{u_n}_x\eqno(6.2)$$
are continuous as functions
from $\R$ into $\L^1_\per\,$, and the identity holds
at every time $t\in\R$, with the exception of the
isolated times where a peakon interaction occurs.

At any time $t$ where no peakon interaction occur in the solution
$u_n$, we define
$\mu^{(n)}_t$ to be the measure with density
$u_n^2(t,\cdot)+{1\over 2}u_{n,x}(t,\cdot)$ w.r.t.~Lebesgue measure.
By Lemmas 5 and 4,
the map $t\mapsto \mu^{(n)}_t$ can be extended by weak continuity
to all times $t\in\R$.
We can now redefine
$$P^{u_n}(t,x)\doteq \int {1\over 2}e^{-|x-y|}\,d\mu_t^{(n)}(y)\,,\qquad\qquad
P^u(t,x)\doteq \int {1\over 2}e^{-|x-y|}\,d\mu_t(y)\,.\eqno(6.3)$$
where $\mu_t$ is the weak limit of the measures $\mu^{(n)}_t$.
Because of the convergence $J\big(u_n(t),\,u(t)\big)\to 0$,
by Lemma 4 the map $t\mapsto \mu_t$ is well defined
and continuous w.r.t.~the weak topology of measures.
Using again
Lemma 4, we can take the limit of (6.2) as $n\to\infty$,
and obtain the identity (1.14), for every $t\in \R$
and $P=P^u$ defined by (6.3).

For each $n$, the total energy
$\mu^{(n)}_t\big(]0,1]\big)=E^{\bar u_n}$
is constant in time and converges to $E^{\bar u}$ as $n\to\infty$.
Therefore we also have
$$\mu_t\big(]0,1]\big)=E^{\bar u}\doteq \int_0^1\big[ \bar u^2(x)+
\bar u_x^2(x)\big]\,dx\qquad\qquad t\in\R\,.\eqno(6.4)$$

To complete the proof of Theorem 1, it now only remains to prove
that the measure $\mu_t$ is absolutely continuous with density
$u^2(t,\cdot)+{1\over 2} u^2_x(t,\cdot)$ w.r.t.~Lebesgue measure,
for a.e.~time $t\in\R$.

In this direction, we recall that, by (1.5), the function $w\doteq
u^2_{n,x}$
satisfies the linear transport equation with source
$$w_t+(uw)_x=(u_n^2-P^{u_n})u_{n,x}\,.\eqno(6.5)$$
Moreover, along any characteristic curve
$t\mapsto\xi(t)$
by (1.3) one has
$${d\over dt}\Big[2\arctan u_{n,x}\big(t,\xi(t,x)\big)\Big]
~=~{2\over 1+u_{n,x}^2}\,
\left[ u_n^2-{u_{n,x}^2\over 2}-P^{u_n}\right]~\leq~-{1\over 2}\,,
\eqno(6.6)$$
whenever $u_{n,x}^2$ is sufficiently large.
For $\ve>0$ small, consider the piecewise affine, $2\pi$-periodic function
(see fig.~2)
$$\vp(\theta)=\left\{\eqalign{ \theta\qquad&\hbox{if}\quad 0\leq\theta\leq 1\,
\cr
1\qquad &\hbox{if}\quad1\leq\theta\leq \pi-\ve\,,\cr
(\pi-\theta)/\ve
\qquad&\hbox{if}\quad \pi-\ve\leq\theta\leq\pi+\ve\,,\cr
-1\qquad&\hbox{if}\quad \pi+\ve\leq\theta\leq 2\pi-1\,,\cr
\theta-2\pi\qquad&\hbox{if}\quad 2\pi-1\leq\theta\leq 2\pi\,.\cr}
\right.
$$
\midinsert
\vskip 10pt
\centerline{\hbox{\psfig{figure=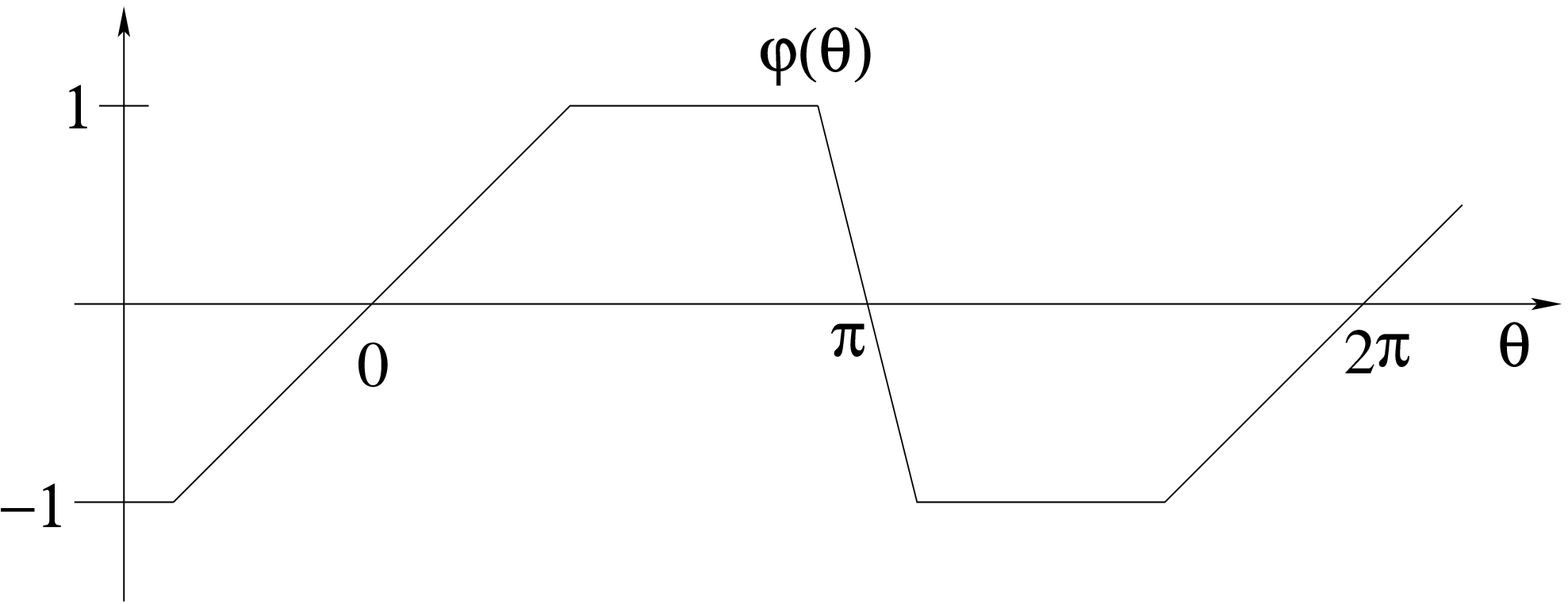,width=10cm}}}
\centerline{\hbox{figure 2}}
\vskip 10pt
\endinsert
\n and define
$$\beta_n(t)=\int_0^1 \vp\big( 2\arctan u_{n,x}(t,x)\big)\, u_{n,x}^2(t,x)
\,dx\,.$$
By (1.3) and (6.6) we now have
$${d\over dt}\beta_n(t)\geq {1\over 4\ve}\int_{\{2\arctan u_{n,x}
\in [\pi-\ve,
\pi]\cup [-\pi, -\pi+\ve]\}} u_{n,x}^2(t,x)\,dx -C\cdot \int_{\{2\arctan u_{n,x}
\in[-1,1]\}}u_{n,x}^2\,dx\eqno(6.7)$$
for some constant $C$, independent of $\ve,n$.
Since all functions $\beta_n$ remain uniformly bounded, by (6.7) for any
time interval $[\tau,\tau']$
we obtain
$$\int_\tau^{\tau'}\int_{\{2\arctan u_{n,x}
\in [\pi-\ve,
\pi]\cup [-\pi, -\pi+\ve]\}} u_{n,x}^2(t,x)\,dx dt\leq
\ve C'\cdot \big(1+\tau'-\tau)\,,\eqno(6.8)$$
where the constant $C'$ depends only on the $H^1_\per$
norm of the functions $u_n$, hence is uniformly valid for all
$n,\ve$.
Because of (6.8), the sequence of functions
$u_{n,x}^2$
is equi-integrable on any domain of the form $[\tau,\tau']\times [0,1]$.
Namely
$$\lim_{\kappa\to\infty}\int_\tau^{\tau'}
\int_{\{x\in[0,1],~ u^2_{n,x} >\kappa\}}
u_{n,x}^2(t,x)\,dxdt=0\,,\eqno(6.9)$$
uniformly w.r.t.~$n$.
By Lemma 4 we already know that $\big\|u_{n,x}^p(t)-u_x^p(t)
\big\|_{\L^1_\per}\to 0$
for every fixed time $t$ and $1\leq p <2$.
Thanks to the equi-integrability condition (6.9) we now have
$$u^2_{n,x}\to u_x^2\qquad\quad \hbox{in}\quad\L^1\big( [\tau,\tau']
\times [0,1]\big)\,.$$
By Fubini's theorem, this implies
$$\lim_{n\to\infty}\int_0^1 u_{n,x}^2(t,x)\,dx=\int_0^1 u_x^2(t,x)\,dx
$$
for a.e.~$t\in [\tau,\tau']$. At every  such time $t$, the measure $\mu_t$
is absolutely continuous and the definition (6.3) coincides with (1.2).
This completes the proof of Theorem 1.
\endproof
\vsk
\n{\medbf 7 - Uniqueness}
\v
Before proving Theorem 3, we remark that the solution satisfying
all conditions in Theorem 1 need not be unique.

\midinsert
\vskip 10pt
\centerline{\hbox{\psfig{figure=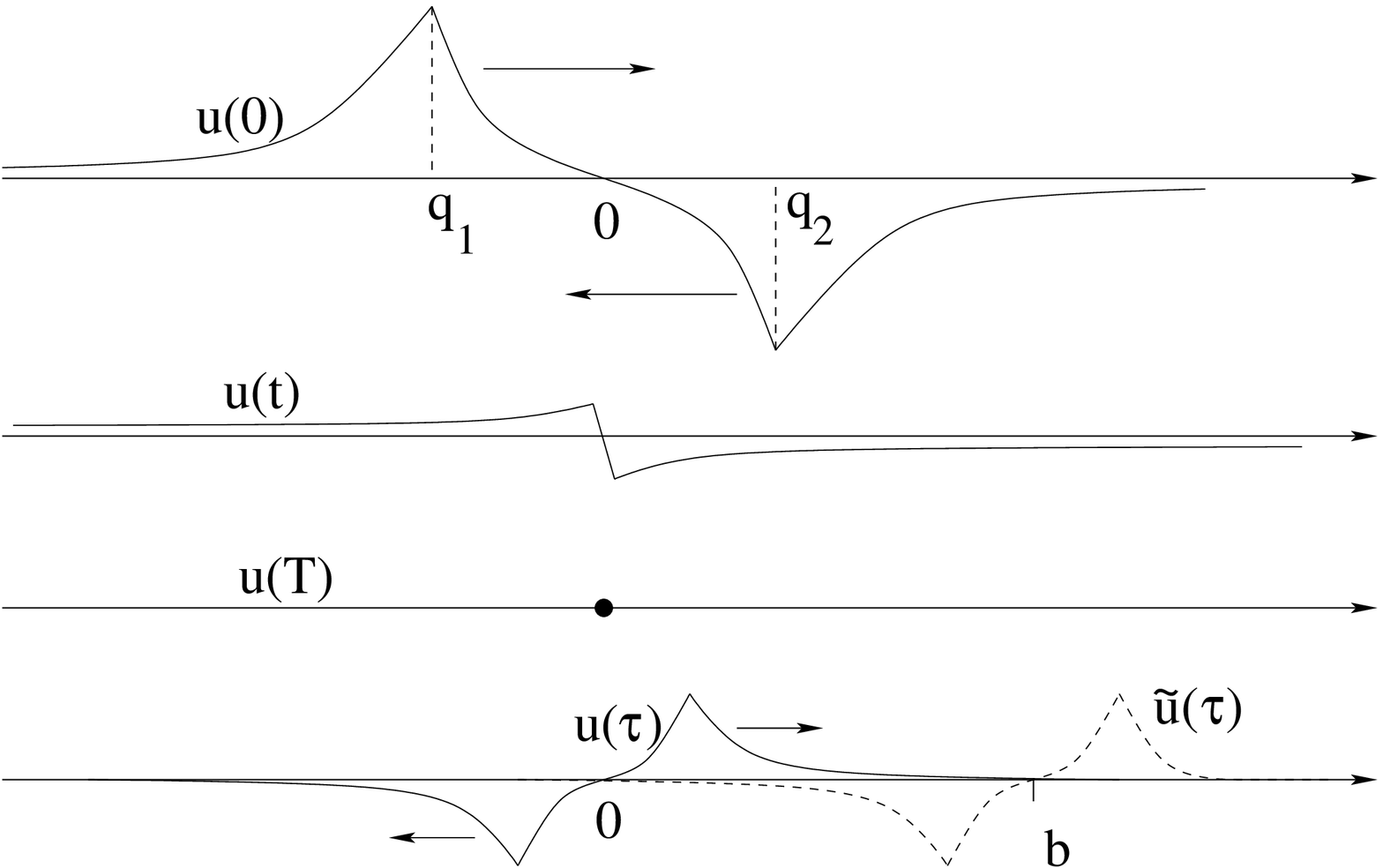,width=12cm}}}
\centerline{\hbox{figure 3}}
\vskip 10pt
\endinsert
\v
\n{\bf Example.} Let $u=u(t,x)$ be a solution containing exactly
two peakons of opposite strengths $p_1(t)=-p_2(t)$, located
at points $q_1(t)=-q_2(t)$ (see fig.~3).  We assume that initially
$p_1(0)>0\,q_1(0)<0$.  At a finite time $T>0$, the two peakons
interact at the origin.  In particular, as $t\to T-$ there holds
$$p_1(t)\to\infty,\qquad p_2(t)\to -\infty,\qquad  q_1(t)\to 0,\qquad
q_2(t)\to 0\,.$$
Moreover, $\big\|u(t)\big\|_{\L^\infty}\to 0$, while
the measure $\mu_t$ approaches a Dirac mass at the origin.
We now have various ways to extend the solution
beyond time $T$:
$$\tilde u(\tau,x)\doteq 0\,,
\eqno(7.1)$$
$$u(\tau,x)=-u(\tau-T,\, -x)\,,\eqno(7.2)$$
Clearly, $\tilde u$ dissipates all the energy, and does
not satisfy the identity (1.15).  The function $u$
in (7.2) is the one constructed by our algorithm in Section 2.
However, there are infinitely many other solutions that
still satisfy (1.15), for example
$$\tilde u(\tau,x)=u(\tau, x-b)\eqno(7.3)$$
where $u$ is as in (7.2) and $b\not= 0$.
The additional condition
in Theorem 3 rules out all of them, because
as $\tau\to T+$, the corresponding measures $\tilde \mu_\tau$
approach a Dirac mass at the point $x=b$, not at the origin.
\vs
We can now give a proof of Theorem 3.
As a first step, we extend our distance $J$
to a larger domain $\D$, consisting of
couples $(u,\mu)$, where $u\in H^1_\per$ and $\mu$ is a
positive (spatially periodic) measure whose absolutely
continuous part has density $u^2+u_x^2$ w.r.t.~Lebesgue measure.
This extension is achieved by continuity:
$$J\big( (u,\mu),\, (\tilde u,\tilde \mu)\big)\doteq
\liminf_{n\to\infty} J(u_n,\tilde u_n)\eqno(7.4)$$
where the infimum is taken over all couple of sequences
$(u_n,\tilde u_n)_{n\geq 1}$ such that
$$\|u_n-u\|_{\L^\infty}\to 0\,,\qquad\qquad
\|\tilde u_n-\tilde u\|_{\L^\infty}\to 0\,,$$
$$u^2_{n,x}\wto \mu\,,\qquad\qquad \tilde u^2_{n,x}\wto\tilde \mu\,.$$
We observe that the flow $\Phi$ constructed in Theorem 2 can be
continuously extended to a locally Lipschitz continuous group
of transformations on the domain $\D$.

Now let $t\mapsto \tilde u(t)$ be a solution of the Cauchy problem
(1.1), (1.13), satisfying all the required conditions.
In particular,  the map $t\mapsto \big(\tilde u(t),\tilde \mu_t\big)$
is Lipschitz continuous w.r.t.~the distance $J$, with values
in the domain $\D$.

Calling $t\mapsto \big(\tilde u(t),\tilde \mu_t\big)
\doteq \Phi_t(\bar u, \bar u^2_x)$
the unique solution of the Cauchy problem obtained
as limit of multi-peakon approximations, we need to
show that $\tilde u(t)=u(t)$
for all $t$.  To fix the ideas, let $t>0$.
By the Lipschitz continuity of the flow,
we can use the error estimate
$$J\Big(\big(\tilde u(t),\tilde \mu_t\big)\,,~\big(u(t), \mu_t\big)
\Big)\leq e^{C_2 t} \int_0^t \liminf_{h\to 0}
{1\over h}\cdot J\Big( \big(\tilde u(\tau+h),\tilde \mu_{\tau+h}\big)\,,~
\Phi_h\big(\tilde u(\tau),\tilde \mu_\tau\big)\Big)\,d\tau
\eqno(7.5)$$
For a proof of (7.5), see pp. 25--27 in [B].
The conditions stated in Theorem 1 now imply that, at almost
every time $\tau$,
the measure $\tilde\mu_t$ is absolutely continuous
and the integrand in (7.5) vanishes.  Therefore $\tilde u(t)=u(t)$
for all $t$.
\endproof
\v
We can now prove that, in multi-peakon solutions,
interactions involving exactly two peakons are the only possible ones.
\v
\n{\bf Corollary.} {\it Let $t\mapsto u(t,\dot)$ be a multi-peakon solution
of the form (1.7), which remains regular
on the open interval $]0,T[\,$.
Assume that at time $T>0$ an interaction occurs, say among the first
$k$ peakons, so that
$$\lim_{t\to T-} q_i(t)=\bar q\qquad\qquad i=1,\ldots,k\,.$$
Then $k=2$.}
\v
\n{\bf Proof.}
We first observe that the Camassa-Holm equations
(1.1) are time-reversible.  In particular, our
proof of Theorem 3 shows that the
solution to a Cauchy problem is unique both forward and backward in time.

Now consider the data
$\big(u(T),\,\mu_T\big)\in\D$, where
$\mu_T$ is the weak limit of the measures $\mu_t$ having density
$u^2(t)+u^2_x(t)$ w.r.t.~Lebesgue measure, as $t\to T-$.
By the analysis in Section 2, we can construct a backward solution
of this Cauchy problem
in terms of exactly two incoming peakons.
By uniqueness, this must coincide with the given solution $u(\cdot)$
for all $t\in [0,T]$.
\endproof
\vsk
\c{\medbf References}
\v
\i{[B]} A.~Bressan,
{\it Hyperbolic Systems of Conservation Laws. The One Dimensional
Cauchy Problem}, Oxford University Press, 2000.
\v
\item{[BC1]} A.~Bressan and A.~Constantin,
Global solutions to the Hunter-Saxton equations,
{\it SIAM J. Math. Anal.}, to appear.
\v
\item{[BC2]} A.~Bressan and A.~Constantin,
Global solutions of the Camassa-Holm equations, submitted.
\v
\item{[BZZ]}
A.~Bressan, P.~Zhang, and Y. Zheng,
On asymptotic variational wave equations,
{\it Arch. Rat. Mech. Anal.}, to appear.
\v
\item{[CH]}
R.~Camassa and D.~D.~Holm,
An integrable shallow water equation with peaked solitons,
{\it Phys. Rev. Lett.}
{\bf 71} (1993), 1661-1664.
\v
\item{[CHK1]}
G.~M.~Coclite, H.~Holden and K.~H.~Karlsen,
Well-posedness for a parabolic-elliptic system,
{\it Discrete Contin. Dynam. Systems.}, to appear.
\v
\item{[CHK2]}
G.~M.~Coclite, H.~Holden and K.~H.~Karlsen,
Global weak solutions to a generalized
hyperelastic-rod wave equation, submitted.
\v
\item{[C1]}
A.~Constantin,
Existence of permanent and breaking waves for a shallow water
equation: a geometric approach, {\it Ann. Inst. Fourier (Grenoble)}
{\bf 50} (2000), 321-362.
\v
\item{[C2]}
A.~Constantin, On the scattering problem for the Camassa-Holm equation,
{\it Proc. Roy. Soc. London Ser. A} {\bf 457} (2001), 953-970
\v
\item{[CE1]}
A.~Constantin and J.~Escher,
Global existence and blow-up for a shallow water equation, {\it
Ann. Scuola Norm. Sup. Pisa} {\bf 26} (1998), 303-328.
\v
\item{[CE2]}
A.~Constantin and J.~Escher, Wave breaking for nonlinear nonlocal
shallow water equations, {\it Acta Mathematica} {\bf 181} (1998), 229-243.
\v
\item{[CM1]}
A.~Constantin and H.~P.~McKean,
A shallow water equation on the circle,
{\it Comm. Pure Appl. Math.}  {\bf 52} (1999), 949-982.
\v
\item{[CM2]}
A.~Constantin and L.~Molinet, Global weak solutions for a shallow
water equation, {\it Comm. Math. Phys.} {\bf 211} (2000), 45-61.
\v
\item{[EG]}
L.~C.~Evans and R.~F.~Gariepy, {\it Measure Theory and Fine
Properties of Functions}, Studies in Advanced Mathematics,
CRC Press, Boca Raton, FL, 1992.
\v
\item{[J]}
R.~S.~Johnson, Camassa-Holm,
Korteweg-de Vries and related models for water waves,
{\it J. Fluid Mech.} {\bf 455} (2002), 63-82.
\v
\i{[V]} C. Villani, {\it Topics in
Optimal Transportation}, Amer. Math. Soc., Providence 2003.
\v
\item{[XZ1}
Z.~Xin and P.~Zhang, On the weak solutions to a shallow water
equation, {\it Comm. Pure Appl. Math.} {\bf 53}  (2000), 1411-1433.
\v
\item{[XZ2]}
Z.~Xin and P.~Zhang, On the uniqueness and large time behavior
of the weak solutions to a shallow water equation, {\it
Comm. Partial Differential Equations}  {\bf 27}  (2002), 1815-1844.

\bye